\documentclass[12pt,a4paper]{article}

\usepackage{amsmath}    
\usepackage{amssymb}
\usepackage{graphicx}   
\usepackage{verbatim}   
\usepackage{color}      

\usepackage{hyperref}   

\usepackage{enumerate}
\usepackage{mathrsfs}
\usepackage{empheq}
\usepackage{bbold}

\usepackage[final]{showlabels}

\usepackage{amsthm}

\newtheorem{lemma}{Lemma}[section]

\newtheorem{theorem}{Theorem}[section]

\renewcommand{\theta}{\vartheta}
\renewcommand{\phi}{\varphi}

\renewcommand{\title}{A Well-Balanced Method for an Unstaggered Central Scheme, the two-space Dimensional Case}

\newcommand{\authorOne}{Yu-Chen Cheng\footnote{Institute of Mathematics, University of Wuerzburg, Wuerzburg, Germany, yu-chen.cheng@stud-mail.uni-wuerzburg.de}}
\newcommand{\authorTwo}{Christian Klingenberg\footnote{Institute of Mathematics, University of Wuerzburg, Wuerzburg, Germany, christian.klingenberg@uni-wuerzburg.de}}
\newcommand{\authorThree}{Rony Touma\footnote{Institute of Computer Science and Mathematics, Lebanese American University, Byblos, Lebanon, rony.touma@lau.edu.lb}}

\begin{document}

\begin{center} \Large
\title

\vspace{1cm}

\date{}
\normalsize

\authorOne, \authorTwo, \authorThree
\end{center}

\begin{abstract}

We develop a second-order accurate central scheme for the two-dimensional hyperbolic system of in-homogeneous conservation laws. The main idea behind the scheme is that we combine the well-balanced {\sl deviation method} with the {\sl Kurganov-Tadmor} (KT) scheme. The approach satisfies the well-balanced property and retains the advantages of KT scheme: Riemann-solver-free and the avoidance of oversampling on the regions between Riemann-fans.
The scheme is implemented and applied to a number of numerical experiments for the Euler equations with gravitational source term and the results are non-oscillatory. Based on the same idea, we construct a semi-discrete scheme where we combine the above two methods and illustrate the maximum principle.



\end{abstract}

\section{Introduction}

The numerical methods  for conservation laws is widely studied in recent years. In 1990, Nessyahu and Tadmor (NT) proposed the central scheme in \cite{ref5} which focuses on the average between two staggered grids. The feature of the NT scheme is to approximate the solutions by integrating over the Riemann fans to avoid solve the Riemann problem at interfaces. After that, authors in \cite{Jiang} developed an unstaggered version of the NT scheme. A modification of the NT scheme was introduced in \cite{ref3} (KT scheme). Compared to the NT scheme which approximate the solution over Riemann fans by integrating on the full-cell size, the KT scheme considers a narrower interval over Riemann fans by the use of the maximal local wave speeds. The idea of the one-dimensional KT scheme is to split the original cell into two intervals. One is the unsmooth region over the whole Riemann fans, and the other is the smooth region between Riemann fans in the original interval. By this approach, the solution is not oversampled at smooth region and can have higher precision. The concept of two dimensional (2-d) extension of the KT scheme is mentioned in \cite{ref8}. Later, a modified 2-d KT scheme was introduced in \cite{ref7}, which considers the maximal local wave speed on the each side of x and y direction. The benefit of this scheme is that the more precise information is used for the solution of unsmooth region.

The set-up of the 2-d KT scheme can be separated into two approaches: one considers the quadrilateral subdomains over the Riemann fans; the other considers the rectangular subdomains. Using rectangular subdomains is easier to formulize and compute for the construction of semi-discrete scheme. There is a list of the 2-d semi-discrete type KT scheme by adopting the rectangular subdomains from the original 2-d KT schemes: \cite{KT2}, \cite{ref8}, \cite{ref9}, and \cite{ref16}. 

In this paper, we focus on in-homogeneous conservation laws, which is given by 
\begin{equation} \label{eq1}
    \partial_t q(x,t)+\nabla_x f(q(x,t))= S(q(x,t)),
\end{equation}
where $q(x,t)=(q_1(x,t), q_2(x,t),...,q_N(x,t))^T$ is an N-vector of conserved quantities in the d-spatial variables $x=(x_1, x_2,..., x_d)$, and $f(q)=(f^1,f^2,...,f^d)$ is a nonlinear flux. $S(q)=(s^1,s^2,...,s^d)$ is a source term. 
Eq.\eqref{eq1} is also known as the balance laws. To solve this equation, a suitable discretization of the source term is needed for construction. There are some attempts for the so-called well-balanced scheme: \cite{ref4}, \cite{ref10}, \cite{ref14}, \cite{hydro}, \cite{WB_BKL}, \cite{upwind2}, \cite{rec_wb2}, \cite{relaxation3_tbc}, \cite{wetdry}, \cite{AMC}, \cite{WenoVol}. 

Inspired by the modified KT scheme and the so-called Deviation method in \cite{ref14}, we construct a well-balanced scheme for the Eq.\eqref{eq1} with a linear source term in next section. In section 3, considering the rectangle subdomains, we derive a semi-discrete scheme by combining the KT scheme and the Deviation method, and we also show the non-oscillation property of this semi-discrete scheme. A number of numerical experiments has been tested by our fully-discrete scheme in section 4. Finally, we end with a conclusion in section 5.  
\section{A new two-dimensional scheme for in-homogeneous conservation laws}\label{sec:scheme}

In this section, we introduce a new two-dimensional scheme by using the combination of the Deviation method and the KT-type scheme. As in \cite{ref4} we begin by deriving the Deviation method at the continuum level for 2-d balance laws with gravity. Then, at the discrete level, we use the the 2-d KT-type scheme in \cite{ref7} to construct a new well-balanced two-dimensional fully-discrete scheme. 
\subsection{Framework of the two-dimensional Deviation method} 
Consider the 2-d balance laws
\begin{equation} \label{4.1}
    \left \{
    \begin{aligned}
    & q_t +f(q)_x + g(q)_y=  S(q,x,y), \;\;\;  (x,y)\in \Omega \subset \mathbb{R}^2, t>0 \\
    & q(x,y,0) = q_0(x,y), \\
    \end{aligned}
    \right.
\end{equation}
where $f(q)$ and $g(q)$ are the fluxes in $x$- and $y$ directions and $S$ is the source term.

Assume $\tilde q$ is a given stationary solution of \eqref{4.1}. Then it satisfies
\begin{equation} \label{4.2}
    f(\tilde q)_x +g(\tilde q)_y = S(\tilde q, x, y).
\end{equation}
Define the deviation $\Delta q = q-\tilde q$. Applying $q= \Delta q+\tilde q$ to \eqref{4.1}, we obtain 
\begin{equation} \label{4.3}
    (\Delta q+\tilde q)_t + f(\Delta q+\tilde q)_x + g(\Delta q+\tilde q)_y = S(\Delta q+ \tilde q).    
\end{equation}
Since $\tilde q$ is a stationary solution, \eqref{4.3} reduces to 
\begin{equation} \label{4.4}
    (\Delta q)_t + f(\Delta q+\tilde q)_x + g(\Delta q+\tilde q)_y = S(\Delta q+ \tilde q).
\end{equation}
Then, we subtract \eqref{4.2} from \eqref{4.4}, and assume that the source term $S(\Delta q+\tilde q)$ in \eqref{4.1}  is a linear functional in terms of the conserved variables, we obtain 
\begin{equation} \label{4.5}
\begin{split}
    (\Delta q)_t + [f(\Delta q+\tilde q) - f(\tilde q)]_x +[g(\Delta q+\tilde q)-g(\tilde q)]_y
      & = S(\Delta q+\tilde q,x,y) - S(\tilde q,x,y) \\
      & = S(\Delta q,x,y).
\end{split}
\end{equation}

\begin{lemma}
\label{lemma1}
Consider the balance law \eqref{4.1} and a given hydrostatic solution $\tilde q$. The deviation quantity $\Delta q$ satisfying the modified balance law \eqref{4.5} maintains the same local speeds as those in the original balance system \eqref{4.1}.
\end{lemma}

\begin{proof}
The proof is the extension of 1D case.  We refer reader to \cite{1D}.
\end{proof}

\subsection{Application of the two-dimensional Kurganov-Tadmor-type scheme} \label{sec:2.2}
Now we apply the idea of KT scheme. The derivation can be separated  into three steps: Reconstruction, Evolution, and Projection. The first two steps are loosely based on the setup in \cite{ref7}. 

\subsubsection{Reconstruction} \label{reconstruction}
Consider the control cell $C_{j,k} = [x_{j-\frac{1}{2}}, x_{j+\frac{1}{2}}]\times[y_{k-\frac{1}{2}}, y_{k+\frac{1}{2}}]$ for all $j, k$. To avoid oscillation, we define a piecewise-linear reconstruction $Q$, 
\begin{equation} \label{4.6}
    Q_{j,k}(x,y,t^n) = (\Delta q)^n_{j,k} + (x-x_j)((\Delta q)_x)^n_{j,k} +(y-y_k)((\Delta q)_y)^n_{j,k}, \;\;\; \forall (x,y)\in C_{j,k},
\end{equation}
where $(\Delta q)_x$ and $(\Delta q)_y$ are the $x$- and $y$ derivatives of $\Delta q$, and apply the $MC$-$\theta$ limiter to estimate the numerical gradient of the solution. The $MC$-$\theta$ limiter is defined as follows,
\begin{equation} \label{4.7}
\begin{split}
    & ((\Delta q)_x)^n_{j,k} = \text{minmod} \left( \theta \Delta_x^+ (\Delta q)^n_{j,k},\, \Delta_x^0 (\Delta q)^n_{j,k},\,\theta \Delta_x^- (\Delta q)^n_{j,k} \right), \\
    & ((\Delta q)_y)^n_{j,k} = \text{minmod} \left( \theta \Delta_y^+ (\Delta q)^n_{j,k},\,\Delta_y^0 (\Delta q)^n_{j,k},\,\theta \Delta_y^- (\Delta q)^n_{j,k} \right),\quad  1\leq \theta \leq 2,
\end{split}
\end{equation}
where the differences $\Delta^\pm_x$, $\Delta^0_x$, $\Delta^\pm_y$, and $\Delta^0_y$ are defined as
\begin{equation} \label{4.8}
\begin{split}
    & \Delta^+_x(\cdot)_{j,k}:= \frac{(\cdot)_{j+1,k}-(\cdot)_{j,k}}{\Delta x},
    \Delta^0_x(\cdot)_{j,k}:= \frac{(\cdot)_{j+1,k}-(\cdot)_{j-1,k}}{2\Delta x},
    \Delta^-_x(\cdot)_{j,k}:= \frac{(\cdot)_{j,k}-(\cdot)_{j-1,k}}{\Delta x}, \\
    & \Delta^+_y(\cdot)_{j,k}:= \frac{(\cdot)_{j,k+1}-(\cdot)_{j,k}}{\Delta y},
    \Delta^0_y(\cdot)_{j,k}:= \frac{(\cdot)_{j,k+1}-(\cdot)_{j,k-1}}{2\Delta y},
    \Delta^-_y(\cdot)_{j,k}:= \frac{(\cdot)_{j,k}-(\cdot)_{j,k-1}}{\Delta y}.
\end{split}
\end{equation}
We will adopt this $(MC-\theta)$ limiter to evaluate the values of the slopes $(\Delta q)_x$ and $(\Delta q)_y$ in the numerical experiments in section \ref{sec:numerics}.

\subsubsection{Evolution} \label{evolution}
\begin{figure}[htbp]
\centering
\includegraphics[width=0.9\textwidth]{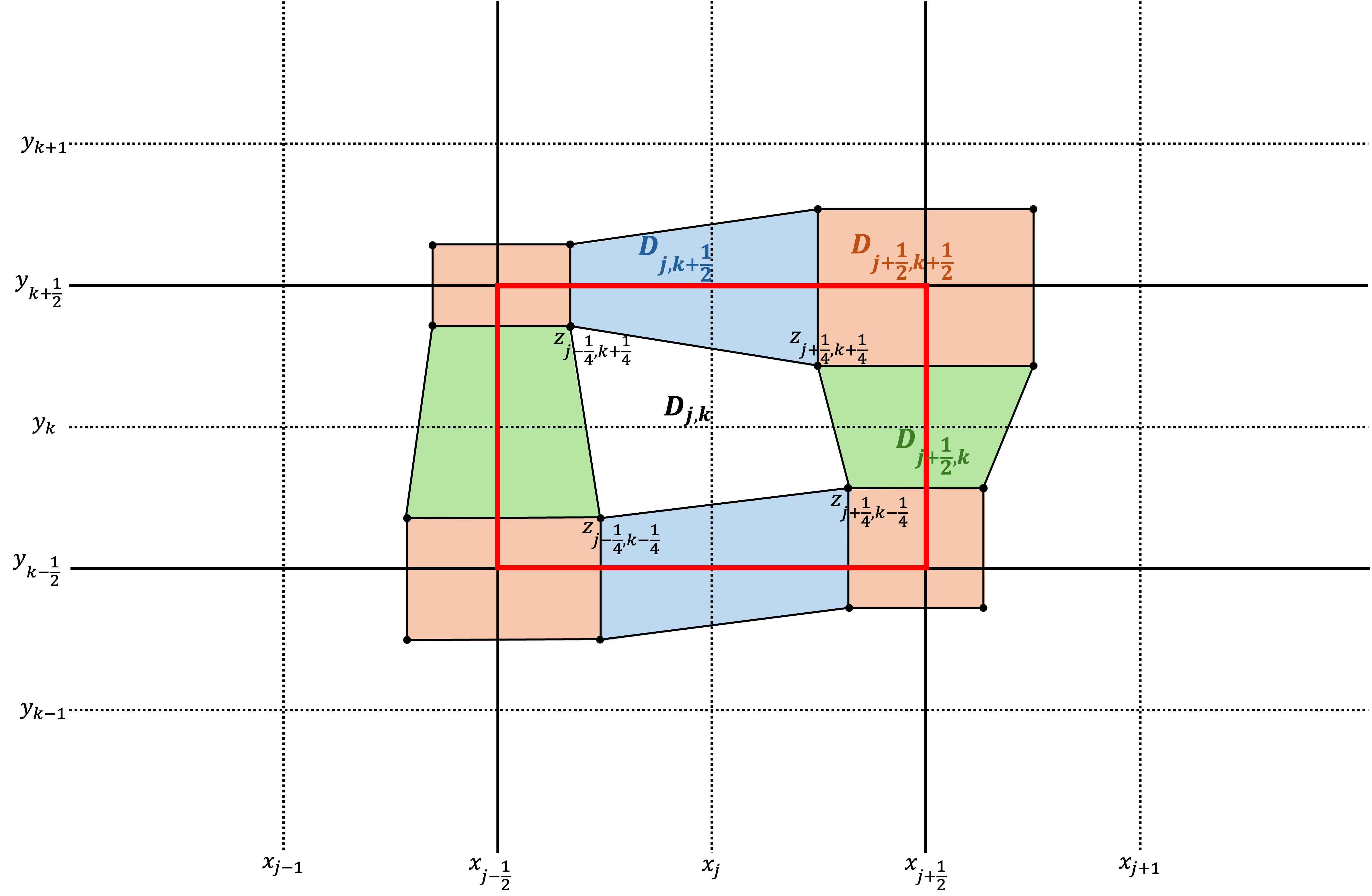}
\caption{Computational cells are split into smooth and unsmooth non-uniform quadrilateral subdomains by the maximal local wave speeds.}
\label{Fig.4.1}
\end{figure}

Before evolving the solution from $t^n$ to the next time step $t^{n+1}$, we need to find the maximal local wave speeds on each side of cell interfaces, which is the most important key of the KT scheme and used to decide the region of the cell over Riemann fans. The local wave speeds at interfaces denoted by $a^\pm_{j+\frac{1}{2},k}$ and $b^\pm_{j,k+\frac{1}{2}}$ in $x$- and $y$ directions, respectively, are determined by the eigenvalues of flux Jacobian:
\begin{equation} \label{4.9}
\begin{split}
    & a^+_{j+\frac{1}{2},k}:= \max \left \{ \lambda_N \left( \frac{\partial}{\partial q} f(q^+_{j+\frac{1}{2},k} ) \right),\; \lambda_N \left( \frac{\partial}{\partial q} f(q^-_{j+\frac{1}{2},k}) \right),\; \epsilon \right\}, \\
    & a^-_{j+\frac{1}{2},k}:= \min \left \{ \lambda_1 \left( \frac{\partial}{\partial q} f(q^+_{j+\frac{1}{2},k} ) \right),\; \lambda_1 \left( \frac{\partial}{\partial q} f(q^-_{j+\frac{1}{2},k}) \right),\; -\epsilon \right\}, \\   
    & b^+_{j,k+\frac{1}{2}}:= \max \left \{ \lambda_N \left( \frac{\partial}{\partial q} g(q^+_{j,k+\frac{1}{2}} ) \right),\; \lambda_N \left( \frac{\partial}{\partial q} g(q^-_{j,k+\frac{1}{2}}) \right),\; \epsilon \right\},  \\  
    & b^-_{j,k+\frac{1}{2}}:= \min \left \{ \lambda_1 \left( \frac{\partial}{\partial q} g(q^+_{j,k+\frac{1}{2}} ) \right),\; \lambda_1 \left( \frac{\partial}{\partial q} g(q^-_{j,k+\frac{1}{2}}) \right),\; -\epsilon \right\},
\end{split}
\end{equation}
where $\lambda_1 < \lambda_2 < \dots < \lambda_N$ are the eigenvalues of the corresponding Jacobians, and $\epsilon$ is a small positive number and
\begin{equation} \label{4.10}
\begin{split}
    & q^-_{j+\frac{1}{2},k}:= q^n_{j,k} + \frac{\Delta x}{2}(q_x)^n_{j,k},\;\;\;\;
    q^+_{j+\frac{1}{2},k}:= q^n_{j+1,k} - \frac{\Delta x}{2}(q_x)^n_{j+1,k},\\
    & q^-_{j,k+\frac{1}{2}}:= q^n_{j,k} + \frac{\Delta y}{2}(q_y)^n_{j,k},\;\;\;\;
    q^+_{j,k+\frac{1}{2}}:= q^n_{j,k+1} - \frac{\Delta y}{2}(q_y)^n_{j,k+1}.
\end{split}
\end{equation}
These local speeds split the control domain $C_{j,k}$ into nine non-uniform subdomains including the unsmooth side subdomains $D_{j,k+\frac{1}{2}}$ and $D_{j+\frac{1}{2},k}$, unsmooth corner subdomain $D_{j+\frac{1}{2},k+\frac{1}{2}}$, and the smooth central subdomain $D_{j,k}$, (see figure \ref{Fig.4.1}).   

The vertices $z_{j\pm\frac{1}{4}, k\pm\frac{1}{4}}$ of these subdomains are estimated by  
\begin{equation} \label{4.11}
\begin{split}
    & z_{j+\frac{1}{4},k+\frac{1}{4}}:=(x_{j+\frac{1}{2}}+\Delta t^n \min(a^-_{j+\frac{1}{2},k},a^-_{j+\frac{1}{2},k+1}), y_{k+\frac{1}{2}}+\Delta t^n \min(b^-_{j,k+\frac{1}{2}},b^-_{j+1,k+\frac{1}{2}})), \\
    & z_{j-\frac{1}{4},k+\frac{1}{4}}:=(x_{j-\frac{1}{2}}+\Delta t^n \max(a^+_{j-\frac{1}{2},k},a^+_{j-\frac{1}{2},k+1}), y_{k+\frac{1}{2}}+\Delta t^n \min(b^-_{j,k+\frac{1}{2}},b^-_{j-1,k+\frac{1}{2}})), \\
    & z_{j-\frac{1}{4},k-\frac{1}{4}}:=(x_{j-\frac{1}{2}}+\Delta t^n \max(a^+_{j-\frac{1}{2},k},a^+_{j-\frac{1}{2},k-1}), y_{k-\frac{1}{2}}+\Delta t^n \max(b^+_{j,k-\frac{1}{2}},b^+_{j-1,k-\frac{1}{2}})), \\
    & z_{j+\frac{1}{4},k-\frac{1}{4}}:=(x_{j+\frac{1}{2}}+\Delta t^n \min(a^-_{j+\frac{1}{2},k},a^-_{j+\frac{1}{2},k-1}), y_{k-\frac{1}{2}}+\Delta t^n \max(b^+_{j,k-\frac{1}{2}},b^+_{j+1,k-\frac{1}{2}})).
\end{split}
\end{equation}

Now we integrate the modified balance law \eqref{4.5} over nine subdomains, respectively, to estimate the intermediate cell-averages at time $t^{n+1}$. 

Consider a general quadrilateral D. The integration of \eqref{4.5} over $D\times[t^n,t^{n+1}]$, 
\begin{equation} \label{4.12}
\begin{split}
    \int^{t^n+1}_{t^n} \iint_D (\Delta q(x,y,t))_t + [f(\Delta q+\tilde q)-f(\tilde q)]_x & + [g(\Delta q +\tilde q) - g(\tilde q)]_y \,dxdydt \\
    & = \int^{t^n+1}_{t^n} \iint_D S(\Delta q,x,y) \,dxdydt.
\end{split}
\end{equation}
can be written as
\begin{equation} \label{4.13}
\begin{split}
    \iint_D \Delta & q(x,y,t^{n+1}) \,dxdy -  \iint_D \Delta q(x,y,t^n) \,dxdy \\
    & + \int^{t^n+1}_{t^n} \iint_D [f(\Delta q+\tilde q)-f(\tilde q)]_x +[g(\Delta q+\tilde q)-g(\tilde q)]_y \,dxdydt  \\
    & \qquad\qquad\qquad\qquad\qquad\qquad\qquad\qquad\quad = \int^{t^{n+1}}_{t^n} \iint_D S(\Delta q,x,y) \,dxdydt. 
\end{split}
\end{equation}
Moving the second and third integrals on the left-hand-side of \eqref{4.13} to the right-hand-side, and dividing the area $|D|$ in \eqref{4.13}, we obtain 
\begin{equation} \label{4.14}
\begin{split}
    \overline{w}^{n+1}_D
    := & \frac{1}{|D|}\iint_D  \Delta q(x,y,t^{n+1}) \,dxdy  \\
    = & \frac{1}{|D|}\iint_D  \Delta q(x,y,t^n)  \,dxdy  \\ 
    & - \frac{1}{|D|}  \int^{t^n+1}_{t^n} \iint_D [f(\Delta q+\tilde q)-f(\tilde q)]_x + [g(\Delta q+\tilde q)-g(\tilde q)]_y \,dxdydt \\
    &\quad + \frac{1}{|D|}\int^{t^{n+1}}_{t^n}  \iint_D S(\Delta q,x,y) \,dxdydt. 
\end{split}
\end{equation}
where $\overline{w}^{n+1}_D$ denotes the new cell-average over D at $t^{n+1}$. 
Applying the divergence theorem to the flux integration, we rewrite the \eqref{4.14} in the form of
\begin{equation} \label{4.16}
\begin{split}
    \overline{w}^{n+1}_D = & \overline{\Delta q}^n_D \\
    & - \frac{1}{|D|}\int^{t^{n+1}}_{t^n}\oint_{\partial D} \eta^x[f(\Delta q + \tilde q)-f(\tilde q)] + \eta^y [g(\Delta q+\tilde q)-g(\tilde q)]\, dxdydt \\ 
     & + \frac{1}{|D|}\int^{t^{n+1}}_{t^n}\iint_D S(\Delta q,x,y)\,dxdydt, 
\end{split}
\end{equation}
where $\eta^x$ and $\eta^y$ are the outward unit normal vectors, and $\overline{\Delta q}^n_D$ define the cell-average over D at  $t^n$.
Then, we discuss the approximations of $\overline{\Delta q}^n_D$, the flux integral, and the source integral on the right-hand-side of \eqref{4.16} separately.

\begin{figure}[htbp]
\centering
\includegraphics[width=0.9\textwidth]{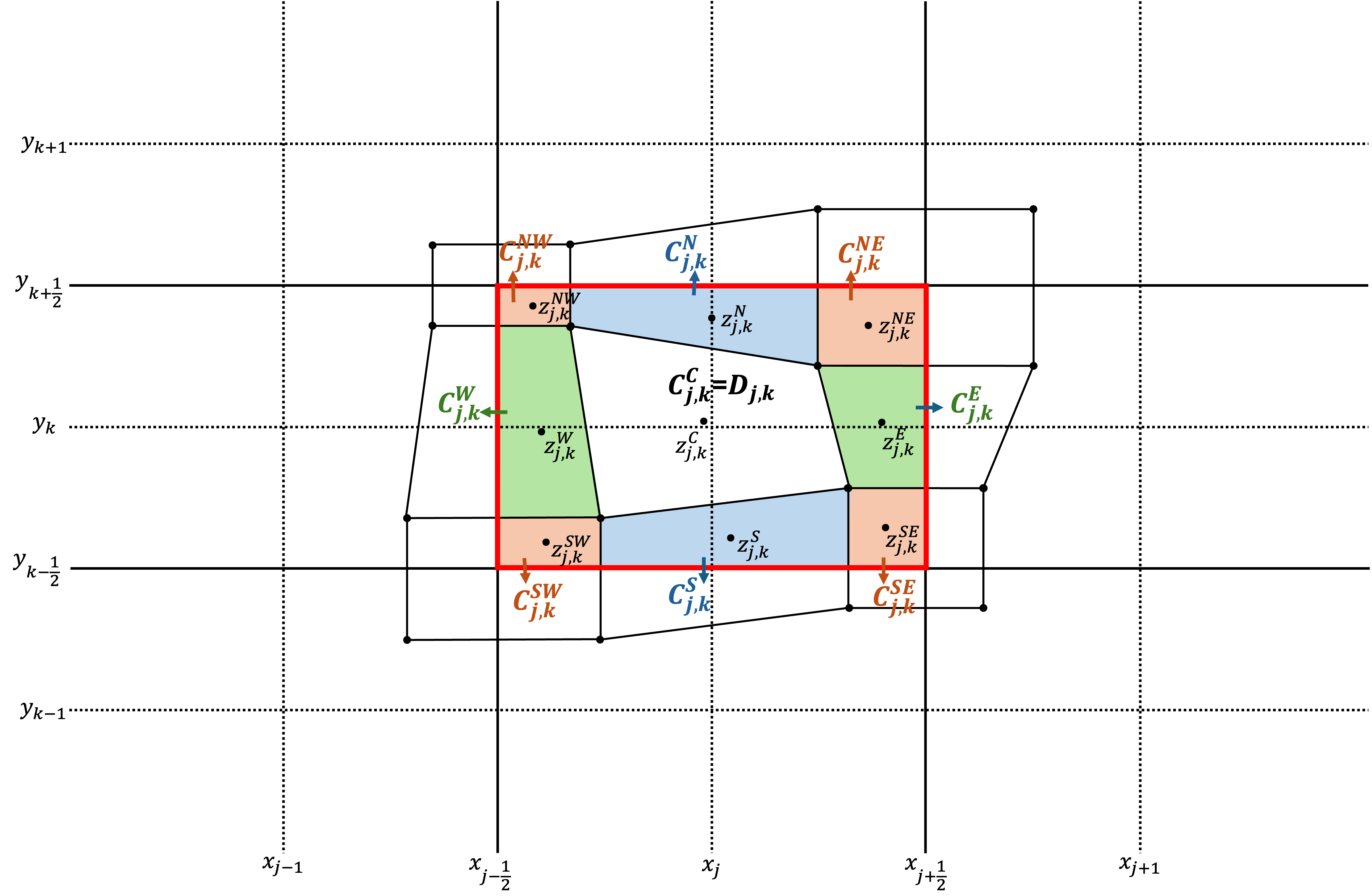}
\caption{Nine subdomains $C^I_{j,k}$ of the original control cell $C_{j,k}$ and their centers of mass.}
\label{Fig.4.3}
\end{figure}

Due to the conservation property of the reconstruction $Q$, the cell-average $\overline{\Delta q}^n_D$ satisfies the following relation
\begin{equation} \label{4.15}
   \overline{\Delta q}^n_D := \frac{1}{|D|}\iint_D \Delta q(x,y,t^n) \,dxdy= \frac{1}{|D|}\iint_D Q(x,y,t^n) \,dxdy.
\end{equation}
We approximate \eqref{4.15} by the following steps:
\begin{enumerate}[1)]
    \item Apply the second order reconstruction \eqref{4.6} to $Q$.
    \item Consider the values over the related subdomains $C^I_{j,k}$,\\ $I=\{E, N, W, S, NE, NW, SE, SW,C\}$, where $C^I_{j,k}$ represent the unions of the subdomain $D$ and the control cell $C_{j,k}$ (see figure \ref{Fig.4.3}).
    \item Use the fact that the cell-average over $C^I_{j,k}$ can be regarded as the values at the centers of mass of $C^I_{j,k}$ denoted by $z^I_{j,k}$ (see figure \ref{Fig.4.3}).
    \item Calculate the weighted averages with the areas of corresponding $C^I_{j,k}$ and their values at the centers of mass, $\Delta q(z^I_{j,k},t)$, and the total area $|D|$ as well.   
\end{enumerate}
The approach here we use for $\overline{\Delta q}^n_D$ is inspired by  \cite{ref7}, we refer the reader to pages A951-A952 for more specific details.  

\begin{figure}[htbp]
\centering
\includegraphics[width=0.9\textwidth]{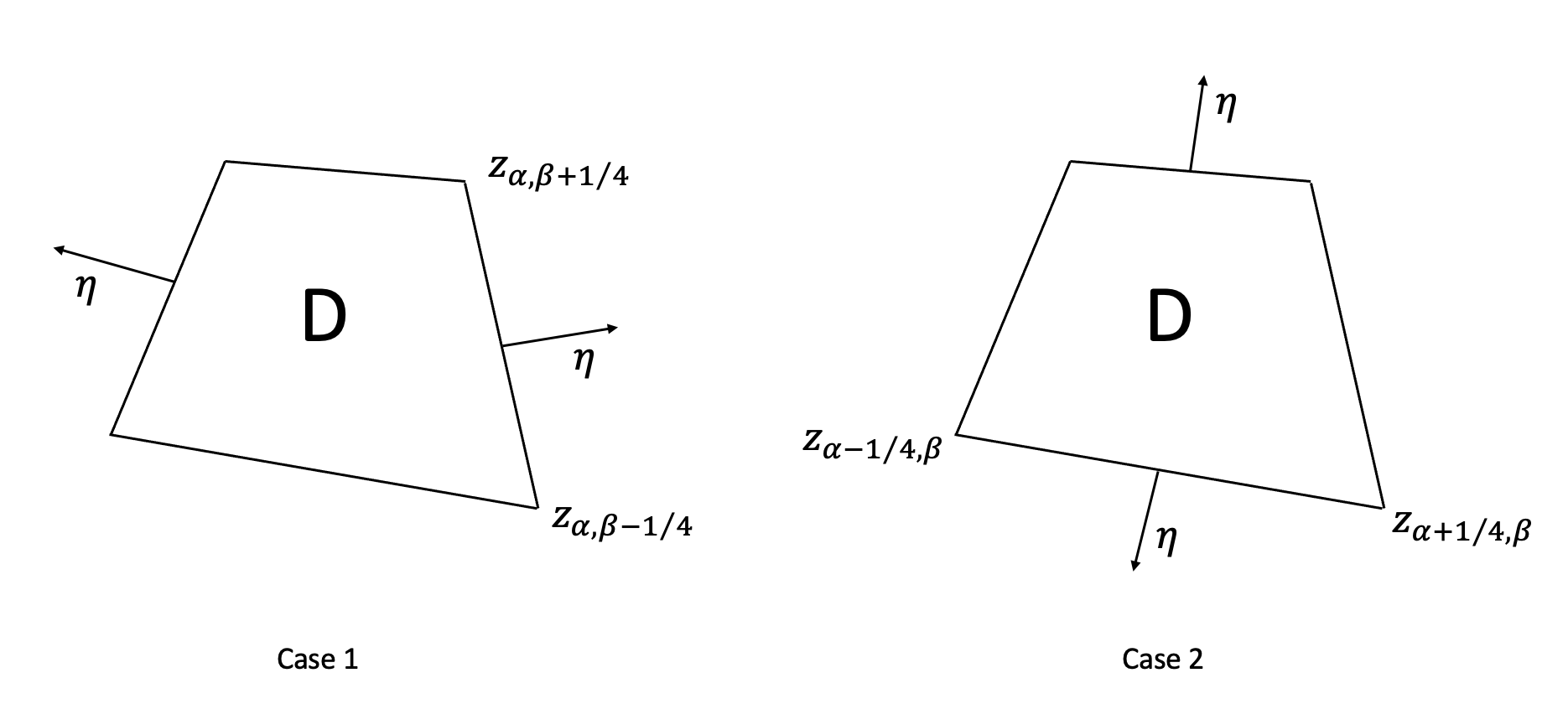}
\caption{A general quadrilateral with normal $\eta$ split into two cases.}
\label{Fig.4.4}
\end{figure}

The evaluation of the flux integral on \eqref{4.16} is also similar to the way in \cite{ref7} that we integrate the fluxes along the four edges of $D$,  and in order to present formula for different directions, we separate them into two cases.\quad As shown in figure \ref{Fig.4.4}, in case 1, we consider the flux integral over the right and left edges, and on the other hand, we consider the flux integral over top and bottom edges in case 2.
\\
\\
\noindent\emph{Case 1.} \\
Consider the edge connecting the nodes $z_{\alpha,\beta-\frac{1}{4}}$ and $z_{\alpha,\beta+\frac{1}{4}}$, where $\alpha = j+\frac{1}{4}$ or $j-\frac{1}{4}$ and $\beta = k$ or $k+\frac{1}{2}$. The flux integration across the edge is defined by 
\begin{equation} \label{4.19}
\begin{split}
    H_{\alpha,\beta} := \frac{1}{\Delta t} \int^{t^{n+1}}_{t^n} \int^{z_{\alpha,\beta+\frac{1}{4}}}_{z_{\alpha,\beta-\frac{1}{4}}}& \Big[ \eta^x_{\alpha,\beta}\left( f(\Delta q+\tilde q)-f(\tilde q) \right) \\
    & +\eta^y_{\alpha,\beta}\left( g(\Delta q+\tilde q)-g(\tilde q) \right) \Big]\, dsdt.
\end{split}
\end{equation}
Applying the midpoint rule in time and the trapezoidal rule in space to \eqref{4.19}, we obtain
\begin{equation} \label{4.20}
\begin{split}
    H_{\alpha,\beta} & = \int^{z_{\alpha,\beta+\frac{1}{4}}}_{z_{\alpha,\beta-\frac{1}{4}}} \left[ \eta^x_{\alpha,\beta}\left( f(\Delta q^{n+\frac{1}{2}}+\tilde q)-f(\tilde q) \right)  
    +\eta^y_{\alpha,\beta}\left( g(\Delta q^{n+\frac{1}{2}}+\tilde q)-g(\tilde q) \right) \right] ds \\
   & = \frac{|z_{\alpha,\beta+\frac{1}{4}}-z_{\alpha,\beta-\frac{1}{4}}|}{2} \left\{ \eta^x_{\alpha,\beta}\left[ f(\Delta q^{n+\frac{1}{2}}_{\alpha,\beta-\frac{1}{4}}+\tilde q_{\alpha,\beta-\frac{1}{4}})-f(\tilde q_{\alpha,\beta-\frac{1}{4}}) \right.\right. \\ 
     &\qquad\qquad\qquad\qquad\qquad\qquad\qquad\qquad \left.+ f(\Delta q^{n+\frac{1}{2}}_{\alpha,\beta+\frac{1}{4}}+\tilde q_{\alpha,\beta+\frac{1}{4}})-f(\tilde q_{\alpha,\beta+\frac{1}{4}}) \right] \\ 
    &\qquad\quad\qquad\qquad\qquad + \eta^y_{\alpha,\beta}\left[ g(\Delta q^{n+\frac{1}{2}}_{\alpha,\beta-\frac{1}{4}}+\tilde q_{\alpha,\beta-\frac{1}{4}})-g(\tilde q_{\alpha,\beta-\frac{1}{4}}) \right. \\ 
    &\qquad\qquad\qquad\qquad\qquad\qquad\qquad\qquad \left. \left. + g(\Delta q^{n+\frac{1}{2}}_{\alpha,\beta+\frac{1}{4}}+\tilde q_{\alpha,\beta+\frac{1}{4}})-g(\tilde q_{\alpha,\beta+\frac{1}{4}}) \right]    \right\},
\end{split}
\end{equation}
where the unit normal vectors $\eta_{\alpha,\beta}$  are
\begin{equation} \label{4.21}
    \eta_{\alpha,k} = \frac{(y_{\alpha,k+\frac{1}{4}}-y_{\alpha,k-\frac{1}{4}},x_{\alpha,k+\frac{1}{4}}-x_{\alpha,k-\frac{1}{4}})}{|z_{\alpha,\beta+\frac{1}{4}}-z_{\alpha,\beta-\frac{1}{4}}|} \qquad
    \text{and} \qquad \eta_{\alpha,k+\frac{1}{2}} =(1,0),
\end{equation}
while $x_{\alpha,\beta}$ and  $y_{\alpha,\beta}$ denote the coordinates of $z_{\alpha,\beta}$.
\\
\\
\noindent\emph{Case 2.} \\
Consider the edge connecting the nodes $z_{\alpha-\frac{1}{4},\beta}$ and $z_{\alpha+\frac{1}{4},\beta}$, where $\alpha = j$ or $j+\frac{1}{2}$ and $\beta = k+\frac{1}{4}$ or $k-\frac{1}{4}$. The flux integration across the edge is defined by 
\begin{equation} \label{4.22}
\begin{split}
    H_{\alpha,\beta} := \frac{1}{\Delta t} \int^{t^{n+1}}_{t^n} \int^{z_{\alpha+\frac{1}{4},\beta}}_{z_{\alpha-\frac{1}{4},\beta}} & \Big[ \eta^x_{\alpha,\beta}\left( f(\Delta q+\tilde q)-f(\tilde q) \right) \\
    & +\eta^y_{\alpha,\beta}\left( g(\Delta q+\tilde q)-g(\tilde q) \right) \Big]\, dsdt.
\end{split}
\end{equation}
Likewise, applying the midpoint rule in time and the trapezoidal rule in space to \eqref{4.22} results in
\begin{equation} \label{4.23}
\begin{split}
    H_{\alpha,\beta} & = \int^{z_{\alpha+\frac{1}{4},\beta}}_{z_{\alpha-\frac{1}{4},\beta}} \left[ \eta^x_{\alpha,\beta}\left( f(\Delta q^{n+\frac{1}{2}}+\tilde q)-f(\tilde q) \right)  
    +\eta^y_{\alpha,\beta}\left( g(\Delta q^{n+\frac{1}{2}}+\tilde q)-g(\tilde q) \right) \right] ds \\
   & = \frac{|z_{\alpha+\frac{1}{4},\beta}-z_{\alpha-\frac{1}{4},\beta}|}{2} \left\{ \eta^x_{\alpha,\beta}\left[ f(\Delta q^{n+\frac{1}{2}}_{\alpha-\frac{1}{4},\beta}+\tilde q_{\alpha-\frac{1}{4},\beta})-f(\tilde q_{\alpha-\frac{1}{4},\beta}) \right.\right. \\ 
     &\qquad\qquad\qquad\qquad\qquad\qquad\qquad\qquad \left.+ f(\Delta q^{n+\frac{1}{2}}_{\alpha+\frac{1}{4},\beta}+\tilde q_{\alpha+\frac{1}{4},\beta})-f(\tilde q_{\alpha+\frac{1}{4},\beta}) \right] \\ 
    &\qquad\quad\qquad\qquad\qquad + \eta^y_{\alpha,\beta}\left[ g(\Delta q^{n+\frac{1}{2}}_{\alpha-\frac{1}{4},\beta}+\tilde q_{\alpha-\frac{1}{4},\beta})-g(\tilde q_{\alpha-\frac{1}{4},\beta}) \right. \\ 
    &\qquad\qquad\qquad\qquad\qquad\qquad\qquad\qquad \left. \left. + g(\Delta q^{n+\frac{1}{2}}_{\alpha+\frac{1}{4},\beta}+\tilde q_{\alpha+\frac{1}{4},\beta})-g(\tilde q_{\alpha+\frac{1}{4},\beta}) \right]    \right\},
\end{split}
\end{equation}
where the unit normal vectors $\eta_{\alpha,\beta}$  are
\begin{equation} \label{4.24}
    \eta_{j,\beta} = \frac{(y_{j-\frac{1}{4},\beta}-y_{j+\frac{1}{4},\beta},x_{j+\frac{1}{4},\beta}-x_{j-\frac{1}{4},\beta})}{|z_{\alpha+\frac{1}{4},\beta}-z_{\alpha-\frac{1}{4},\beta}|} \qquad
    \text{and} \qquad \eta_{j+\frac{1}{2},\beta} =(0,1).
\end{equation}
To obtain the solution at the nodes $z_{j\pm\frac{1}{4},k\pm\frac{1}{4}}$, we apply the Taylor expansion twice (once in time and once in space), and then get the second-order approximation
\begin{equation} \label{4.25}
\begin{split}
    \Delta q^{n+\frac{1}{2}}_{j\pm\frac{1}{4},k\pm\frac{1}{4}} := \Delta q^n_{j\pm\frac{1}{4},k\pm\frac{1}{4}} + & \frac{\Delta t}{2}((\Delta q)_t)^n_{j\pm\frac{1}{4},k\pm\frac{1}{4}} \\
    = \Delta q^n_{j\pm\frac{1}{4},k\pm\frac{1}{4}} - & \frac{\Delta t}{2}\left[ \left[ (f(\Delta q+\tilde q)-f(\tilde q))_x \right]^n_{j\pm\frac{1}{4},k\pm\frac{1}{4}} \right.\\ 
    & \quad\; +\left. \left[ (g(\Delta q+\tilde q)-g(\tilde q))_y \right]^n_{j\pm\frac{1}{4},k\pm\frac{1}{4}}\right] + \frac{\Delta t}{2} S(\Delta q)^n_{j\pm\frac{1}{4},k\pm\frac{1}{4}} \\
    := \Delta q^n_{j\pm\frac{1}{4},k\pm\frac{1}{4}} - & \frac{\Delta t}{2}\left[ \left[ (f(\Delta q+\tilde q)-f(\tilde q))_x \right]^n_{j,k} 
     + \left[ (g(\Delta q+\tilde q)-g(\tilde q))_y \right]^n_{j,k}\right] \\
     & \qquad\qquad\qquad\qquad\qquad\qquad\qquad\quad + \frac{\Delta t}{2} S(\Delta q)^n_{j\pm\frac{1}{4},k\pm\frac{1}{4}}, 
\end{split}
\end{equation}
where the values of $\Delta q^n_{j\pm\frac{1}{4},k\pm\frac{1}{4}}$ are also estimated by the Taylor expansion:
\begin{equation} \label{4.26}
    \Delta q^n_{j\pm\frac{1}{4},k\pm\frac{1}{4}} := \Delta q^n_{j,k} + (x_{j\pm\frac{1}{4},k\pm\frac{1}{4}}-x_j)(\Delta q_x)^n_{j,k} + (y_{j\pm\frac{1}{4},k\pm\frac{1}{4}}-y_k)(\Delta q_y)^n_{j,k},
\end{equation}
and the slopes $(F(\Delta q)_x)^n_{j,k}:=((f(\Delta q+\tilde q)-f(\tilde q))_x)^n_{j,k}$ and $(G(\Delta q)_y)^n_{j,k}:=((g(\Delta q+\tilde q)-g(\tilde q))_y)^n_{j,k}$ which can be obtained by using the $(MC-\theta)$ limiter:
\begin{equation} \label{4.27}
\begin{split}
    (F(\Delta q)_x)^n_{j,k} = \text{minmod} \left( \theta \Delta_x^+ \left(F(\Delta q^n_{j,k})\right),\, \Delta_x^0 \left(F(\Delta q^n_{j,k})\right),\,\theta \Delta_x^- \left(F(\Delta q^n_{j,k})\right) \right),& \\
    (G(\Delta q)_y)^n_{j,k} = \text{minmod} \left( \theta \Delta_y^+ \left(G(\Delta q^n_{j,k})\right),\, \Delta_y^0 \left(G(\Delta q^n_{j,k})\right),\,\theta \Delta_y^- \left(G(\Delta q^n_{j,k})\right) \right),& \\
    \theta \in [1 & ,2]. 
\end{split}
\end{equation}
 By the approximations \eqref{4.20} and \eqref{4.23}, the sum of the flux along the edges yields the result of the flux integral on \eqref{4.16},
\begin{equation} \label{4.28}
\begin{split}
    \frac{1}{|D_{\alpha,\beta}|}\int^{t^{n+1}}_{t^n}\oint_{\partial D_{\alpha,\beta}}  \eta^x[f(\Delta q + & \tilde q)  -f(\tilde q)] + \eta^y [g(\Delta q+\tilde q)-g(\tilde q)] dxdydt \\
    & = \frac{\Delta t^n}{|D_{\alpha,\beta}|} \left[ H_{\alpha+\frac{1}{4},\beta} - H_{\alpha-\frac{1}{4},\beta} +H_{\alpha,\beta+\frac{1}{4}} - H_{\alpha,\beta-\frac{1}{4}} \right],
\end{split}
\end{equation}
where $\alpha = j$ or $j\pm\frac{1}{2}$ and $\beta =k$ or $k\pm\frac{1}{2}$. 

Finally, we consider the last integral in \eqref{4.16}. To evaluate the average of the source term, we use the midpoint rule in time and consider the value at center of mass $\Delta q(z^I_{j,k},t)$ (which is defined in the page A951-A952 of \cite{ref7} and used for approximating $\overline{\Delta q}^n_D$). The approximation of the source term over the central subdomain $D_{j,k}$ is then defined by        
\begin{equation} \label{4.29}
\begin{split}
    \frac{1}{|D_{j,k}|} & \int^{t^{n+1}}_{t^n} \iint_{D_{j,k}} S(\Delta  q(x,y,t))dxdydt \\
    := & \frac{\Delta t}{|D_{j,k}|}\iint_{D_{j,k}} S (\Delta  q(x,y,t^{n+\frac{1}{2}}))dxdy \\
    \;\quad = & \frac{\Delta t}{|D_{j,k}|} 
    |D_{j,k}| S(\Delta q(z^C_{j,k},t^{n+\frac{1}{2}})) 
    =: \Delta t S_{D_{j,k}}(\Delta q),
\end{split}
\end{equation}
The approximations of the side subdomain $D_{j+\frac{1}{2},k}$ and $D_{j,k+\frac{1}{2}}$, and the corner $D_{j+\frac{1}{2},k+\frac{1}{2}}$ are similar to \eqref{4.15}, and the detailed calculations can be found in \ref{append.A}.   

To summarize this evolution step, we substitute the results of the cell-average $\overline{(\Delta q)^n_D}$, the flux integral \eqref{4.28}, and the source integral \eqref{4.29}, \eqref{4.30}-\eqref{4.32} for \eqref{4.16},
\begin{equation} \label{4.34}
\begin{split}
    \overline{w}^{n+1}_{D_{\alpha,\beta}} = \overline{\Delta q^n_{D_{\alpha,\beta}}} -  \frac{\Delta t}{|D_{\alpha,\beta}|}  \left[ H_{\alpha+\frac{1}{4},\beta} - H_{\alpha-\frac{1}{4},\beta} + H_{\alpha,\beta+\frac{1}{4}} - \right.&\left. H_{\alpha,\beta-\frac{1}{4}} \right] \\
       & + \Delta t S_{D_{\alpha,\beta}}(\Delta q),
\end{split}
\end{equation}
where $\alpha=j$ or $j\pm\frac{1}{2}$ and $\beta=k$ or $k\pm\frac{1}{2}$.

\subsubsection{Projection} \label{projection}
At the final procedure, we project the intermediate solutions $\overline{w}^{n+1}_D$ in \eqref{4.34} back onto the original uniform cell $C_{j,k}$. As stated in the reconstruction step above, to smooth the solution, we need to define a piecewise-linear reconstruction $\widetilde W^{n+1}(x,y)$ for $\overline{w}^{n+1}_D$.

We firstly consider the central smooth subdomain $D_{j,k}$. Since $D_{j,k}\subset C_{j,k}$, the intermediate solution $\overline{w}^{n+1}_{D_{j,k}}$ is smooth enough and does not need a reconstruction.\quad Thus, 
\begin{equation} \label{4.36}
    \widetilde W^{n+1}(x,y) = \overline{w}^{n+1}_{D_{j,k}}\;\;\;\; \text{for}\;\; (x,y)\in D_{j,k}.
\end{equation}
\par Next, we consider the unsmooth subdomains $D_{\alpha,\beta}$ with $(\alpha,\beta)=(j+\frac{1}{2},k)$, $ (j+\frac{1}{2},k+\frac{1}{2})$, and $(j,k+\frac{1}{2})$, and define the reconstruction of the intermediate solutions over $D_{\alpha,\beta}$ by
\begin{equation} \label{4.37}
    \widetilde W^{n+1}_{D_{\alpha,\beta}} 
    = \overline{w}^{n+1}_{D_{\alpha,\beta}}
    + (x - z^{n,x}_{D_{\alpha,\beta}})(w_x)^{n+1}_{D_{\alpha,\beta}}
    + (y - z^{n,y}_{D_{\alpha,\beta}})(w_y)^{n+1}_{D_{\alpha,\beta}}
\end{equation}
where $z^{n,x(y)}_{D_{\alpha,\beta}}$ denotes the coordinates of the centers of mass of the domain $D_{\alpha,\beta}$, and the spatial derivatives are determined by 
\begin{equation} \label{4.38}
\begin{split}
    (w_x)^{n+1}_{D_{\alpha,\beta}}
    & = \text{minmod}
    (\theta
    \frac{\overline{w}^{n+1}_{D_{\alpha,\beta}}-\overline{w}^{n+1}_{D_{\alpha-\frac{1}{2},\beta}}}{z^{n,x}_{D_{\alpha,\beta}}-z^{n,x}_{D_{\alpha-\frac{1}{2},\beta}}},
    \frac{\overline{w}^{n+1}_{D_{\alpha+\frac{1}{2},\beta}} - \overline{w}^{n+1}_{D_{\alpha-\frac{1}{2},\beta}}}{z^{n,x}_{D_{\alpha+\frac{1}{2},\beta}}-z^{n,x}_{D_{\alpha-\frac{1}{2},\beta}}}, 
    \theta \frac{\overline{w}^{n+1}_{D_{\alpha+\frac{1}{2},\beta}}-\overline{w}^{n+1}_{D_{\alpha,\beta}}}{z^{n,x}_{D_{\alpha+\frac{1}{2},\beta}}-z^{n,x}_{D_{\alpha,\beta}}}
    ) \\
    (w_y)^{n+1}_{D_{\alpha,\beta}}
    & = \text{minmod}
    (\theta
    \frac{\overline{w}^{n+1}_{D_{\alpha,\beta}}-\overline{w}^{n+1}_{D_{\alpha,\beta-\frac{1}{2}}}}{z^{n,y}_{D_{\alpha,\beta}}-z^{n,y}_{D_{\alpha,\beta-\frac{1}{2}}}},
    \frac{\overline{w}^{n+1}_{D_{\alpha,\beta+\frac{1}{2}}} - \overline{w}^{n+1}_{D_{\alpha,\beta-\frac{1}{2}}}}{z^{n,y}_{D_{\alpha,\beta+\frac{1}{2}}}-z^{n,y}_{D_{\alpha,\beta-\frac{1}{2}}}}, 
    \theta \frac{\overline{w}^{n+1}_{D_{\alpha,\beta+\frac{1}{2}}}-\overline{w}^{n+1}_{D_{\alpha,\beta}}}{z^{n,y}_{D_{\alpha,\beta+\frac{1}{2}}}-z^{n,y}_{D_{\alpha,\beta}}}
    ),\\
    &\qquad\qquad\qquad\qquad\qquad\qquad\qquad\qquad\qquad\qquad\qquad\qquad\qquad 1\leq\theta\leq2.
\end{split}
\end{equation}
The computation of the centers of mass of domain $D_{\alpha, \beta}$ are similar to $z^I_{j,k}$, which can be found in Appendix A of \cite{ref7}. 

Then, with the defined reconstruction $\widetilde W^{n+1}(x,y)$, we obtain the new cell-average over the original cell $C_{j,k}$ at time $t^{n+1}$ by    
\begin{equation}  \label{4.35}
\begin{split}
    (\Delta q)^{n+1}_{j,k} &= \frac{1}{\Delta x \Delta y} \iint_{C_{j,k}}\widetilde W^{n+1}(x,y) \;dxdy \\
    &= \frac{1}{\Delta x\Delta y}
    \bigg[  |D_{j,k}| \overline{w}^{n+1}_{D_{j,k}} + |C^E_{j,k}|\widetilde W^{n+1}_{D_{j+\frac{1}{2},k}} + |C^{NE}_{j,k}|\widetilde W^{n+1}_{D_{j+\frac{1}{2},k+\frac{1}{2}}}  \\
 &\qquad\qquad +  |C^N_{j,k}|\widetilde W^{n+1}_{D_{j,k+\frac{1}{2}}} + |C^{NW}_{j,k}|\widetilde W^{n+1}_{D_{j-\frac{1}{2},k+\frac{1}{2}}} +|C^W_{j,k}|\widetilde W^{n+1}_{D_{j-\frac{1}{2},k}}  \\
&\qquad\qquad +  |C^{SW}_{j,k}|\widetilde W^{n+1}_{D_{j-\frac{1}{2},k-\frac{1}{2}}} +|C^S_{j,k}|\widetilde W^{n+1}_{D_{j,k-\frac{1}{2}}} +|C^{SE}_{j,k}|\widetilde W^{n+1}_{D_{j+\frac{1}{2},k-\frac{1}{2}}} \bigg].
\end{split}
\end{equation}

Finally, adding the hydrostatic solution $\tilde q_{j,k}$ to the computed solution $(\Delta q)^{n+1}_{j,k}$, we get the desired solution $q^{n+1}_{j,k}$,
\begin{equation}
    q^{n+1}_{j,k} = (\Delta q)^{n+1}_{j,k} + \tilde q_{j,k}.
\end{equation}

\section{2D Semi-discrete scheme}\label{sec:semi}

In this section, we aim to construct a semi-discrete scheme which preserves the characteristics of both the Deviation method and the KT scheme. In the evolution step above in section \ref{sec:scheme}, we split the control cell to nine subdomains by using the values of local wave speeds and then approximate the cell-average over each subdomain. Likewise, the wave speeds are again used to decide the ranges of the subdomains in this section. Compared to the set-up in section 2, instead of trapezoid subdomains, we consider the rectangle subdomains in the following construction, which are easier to be formalized and computed. 

In order to construct the semi-discrete scheme, we first follow the three-steps structure in section \ref{sec:2.2} to build the fully-discrete scheme. Then we derive the semi-discrete scheme. In this section we extend ideas from \cite{ref8}, \cite{ref9} and \cite{ref16}.

\subsection{Construction of the fully-discrete scheme}
\subsubsection{Reconstruction} 
In this step, we use the same piecewise linear interpolant as in \eqref{4.6}, and the same approximations of the derivatives in the  $x$- and $y$ directions as in \eqref{4.7}-\eqref{4.8}.
\subsubsection{Evolution} 
We adopt the same approximation in \eqref{4.10} to compute the reconstructed values at the interfaces, and set $\epsilon=0$ in \eqref{4.9} to define the local wave speeds. 

\begin{figure}[htbp]
\centering
\includegraphics[width=0.9\textwidth]{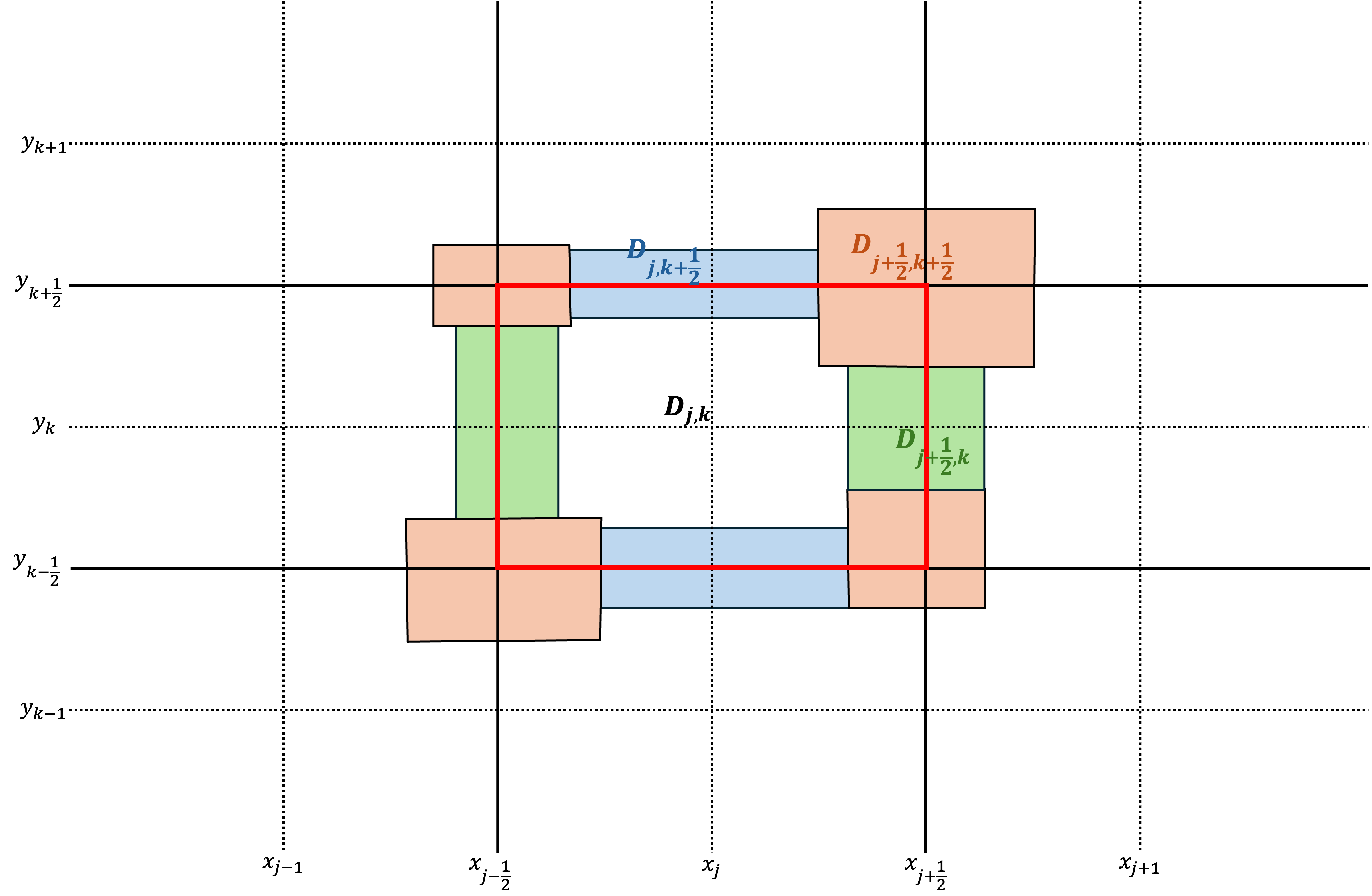}
\caption{Computational cells are split into smooth and unsmooth non-uniform rectangular subdomains by the maximal local wave speeds}
\label{Fig.4.5}
\end{figure}

As we mentioned earlier, to estimate the values over the Riemann fans, we split the control cell to nine non-uniform rectangle subdomains, outlined in figure \ref{Fig.4.5}.

The specific definitions of ranges of the subdomains are given in page 716 of \cite{ref9}. 
Next, we apply \eqref{4.14} to the subdomains $D_{j+\frac{1}{2},k}$-, $D_{j,k+\frac{1}{2}}$-, $D_{j+\frac{1}{2},k+\frac{1}{2}}$- and $D_{j,k}\times [t^n,t^{n+1}]$  to compute the intermediate values at time $t^{n+1}$ denoted by $\overline{w}^{n+1}_{j+\frac{1}{2},k}$, $\overline{w}^{n+1}_{j,k+\frac{1}{2}}$, $\overline{w}^{n+1}_{j+\frac{1}{2},k+\frac{1}{2}}$ and $\overline{w}^{n+1}_{j,k}$ respectively. The resultant integration can be found in \ref{append.B}.

\subsubsection{Projection} 
Finally, we project the intermediate cell-average $\overline{w}^{n+1}_D$ back onto the uniform cell  $C_{j,k}=[x_{j-\frac{1}{2}},x_{j+\frac{1}{2}}]\times[y_{k-\frac{1}{2}},y_{k+\frac{1}{2}}]$ with a suitable reconstruction $\widetilde W^{n+1}(x,y)$ for unsmooth subdomains. For the subdomain $D_{j+\frac{1}{2},k}$, the reconstruction is defined by
\begin{equation} \label{4.50}
\begin{split}
    \widetilde W^{n+1}_{j+\frac{1}{2},k} = \overline{w}^{n+1}_{j+\frac{1}{2},k} 
    + (x-x^{n+1}_{j+\frac{1}{2},k})(w_x)^{n+1}_{j+\frac{1}{2},k}
    + (y-y^{n+1}_{j+\frac{1}{2},k})(w_y)^{n+1}_{j+\frac{1}{2},k},    
\end{split}
\end{equation}
where
\begin{equation}
\begin{split}
    & x^{n+1}_{j+\frac{1}{2},k} := x_{j+\frac{1}{2}} + \frac{a^-_{j+\frac{1}{2},k} + a^+_{j+\frac{1}{2},k}}{2} \Delta t, \\ 
    & y^{n+1}_{j+\frac{1}{2},k} := y_k + \frac{B^+_{j+\frac{1}{2},k-\frac{1}{2}} + B^-_{j+\frac{1}{2},k+\frac{1}{2}}}{2} \Delta t,
\end{split}
\end{equation}
are the coordinates of the center of $D_{j+\frac{1}{2},k}$, and $(w_x)^{n+1}_{j+\frac{1}{2},k}$ and $(w_y)^{n+1}_{j+\frac{1}{2},k}$ are the $x$ and $y$ numerical derivatives. The reconstructions $\widetilde W^{n+1}_{j,k+\frac{1}{2}}$ and $\widetilde W^{n+1}_{j+\frac{1}{2},k+\frac{1}{2}}$ are defined in the similar way. 

Then, the solution $(\Delta q^{n+1}_{j,k})$ reduces to
\begin{equation} \label{4.49}
\begin{split}
    (\Delta q)^{n+1}_{j,k} =& \frac{1}{\Delta x\Delta y} \times \\
    &\times \Bigg( \sum_{E,W}  \iint_{C^{E(W)}_{j,k}} \widetilde W^{n+1}_{j\pm\frac{1}{2},k}(x,y)dxdy
    + \sum_{N,S} \iint_{C^{N(S)}_{j,k}} \widetilde W^{n+1}_{j,k\pm\frac{1}{2}}(x,y)dxdy \\
    &\qquad\qquad\quad + \sum_{NE,NW,SW,SE} \iint_{C^{NE(NW)(SW)(SE)}_{j,k}} \widetilde W^{n+1}_{j\pm\frac{1}{2},k\pm\frac{1}{2}}(x,y) dx dy \\ 
    &\qquad\qquad\qquad\qquad\qquad\qquad\qquad\qquad\qquad + \iint_{D_{j,k}} \overline{w}^{n+1}_{j,k}(x,y)dxdy   \Bigg),
\end{split}
\end{equation}
where subdomains $C^I_{j,k}$, $I =\{ E, W, N, S, NE, NW, SW,SE \}$ denote the union of subdomains $D$ and control cell $C_{j,k}$.

Next, we proceed to the construction of semi-discrete scheme.

\subsection{New semi-discrete scheme}\label{sec:3.2}
With the formula \eqref{4.49}, to obtain the semi-discrete scheme we consider the following derivative, 
\begin{equation} \label{4.52}
\begin{split}
    & \frac{d(\Delta q)_{j,k}(t)}{dt} \\ 
    &\quad= \lim_{\Delta t\to 0}\frac{(\Delta q)^{n+1}_{j,k}-(\Delta q)^n_{j,k}}{\Delta t} \\ 
     &\quad = \lim_{\Delta t\to 0} \frac{1}{\Delta t} \times \\
     &\qquad \times \Bigg[\frac{1}{\Delta x\Delta y} \bigg( \sum_{E,W} \iint_{C^{E(W)}_{j,k}} \widetilde W^{n+1}_{j\pm\frac{1}{2},k}(x,y)dxdy 
     + \sum_{N,S} \iint_{C^{N(S)}_{j,k}} \widetilde W^{n+1}_{j,k\pm\frac{1}{2}}(x,y) dx dy \\
     &\qquad\qquad\qquad\qquad\qquad + \sum_{NE,NW,SE,SW} \iint_{C^{NE(NW)(SE)(SW)}_{j,k}} \widetilde W^{n+1}_{j\pm\frac{1}{2},k\pm\frac{1}{2}}(x,y) dx dy \\   &\qquad\qquad\qquad\qquad\qquad\quad\qquad\qquad\qquad\qquad + \iint_{D_{j,k}} \overline{w}^{n+1}_{j,k}(x,y) dx dy   \bigg) -\Delta q^n \Bigg].
\end{split}
\end{equation}
By using the conservation property and some technical derivations (shown in \ref{append.C}), the semi-discrete scheme then take the form that 
\begin{equation} \label{4.65}
    \frac{d}{dt} (\Delta q)_{j,k}(t) = 
    - \frac{H^x_{j+\frac{1}{2},k}(t)-H^x_{j-\frac{1}{2},k}(t)}{\Delta x}
    - \frac{H^y_{j,k+\frac{1}{2}}(t)-H^y_{j,k-\frac{1}{2}}(t)}{\Delta y}
    + S((\Delta q)_{j,k}(t)),
\end{equation}
with the numerical fluxes
\begin{equation} \label{4.66}
\begin{split}
    H^x_{j+\frac{1}{2},k}  = & \frac{a^+_{j+\frac{1}{2},k}F((\Delta q)^E_{j,k})-a^-_{j+\frac{1}{2},k}F((\Delta q)^W_{j+1,k})}{a^+_{j+\frac{1}{2},k}-a^-_{j+\frac{1}{2},k}}  \\
    &\qquad\qquad\qquad\qquad\qquad + \frac{a^+_{j+\frac{1}{2},k}a^-_{j+\frac{1}{2},k}}{a^+_{j+\frac{1}{2},k}-a^-_{j+\frac{1}{2},k}} \left[(\Delta q)^W_{j+1,k}-(\Delta q)^E_{j,k} \right], \\
    H^y_{j,k+\frac{1}{2}} = & \frac{b^+_{j,k+\frac{1}{2}}G((\Delta q)^N_{j,k})-b^-_{j,k+\frac{1}{2}}G((\Delta q)^S_{j,k+1})}{b^+_{j,k+\frac{1}{2}}-b^-_{j,k+\frac{1}{2}}} \\
    &\qquad\qquad\qquad\qquad\qquad + \frac{b^+_{j,k+\frac{1}{2}}b^-_{j,k+\frac{1}{2}}}{b^+_{j,k+\frac{1}{2}}-b^-_{j,k+\frac{1}{2}}} \left[(\Delta q)^S_{j,k+1}-(\Delta q)^N_{j,k} \right],
\end{split}
\end{equation}
where
\begin{equation} \label{4.60}
\begin{split}
    F((\Delta q)^W_{j+1,k}) &= f((\Delta q)^W_{j+1,k}+\tilde q_{j+\frac{1}{2},k})-f(\tilde q_{j+\frac{1}{2},k}), \\
    F((\Delta q)^E_{j,k}) &= f((\Delta q)^E_{j,k}+\tilde q_{j+\frac{1}{2},k})-f(\tilde q_{j+\frac{1}{2},k}), \\ 
    G((\Delta q)^S_{j,k+1}) &= g((\Delta q)^S_{j,k+1}+\tilde q_{j,k+\frac{1}{2}})-g(\tilde q_{j,k+\frac{1}{2}}), \\
    G((\Delta q)^N_{j,k}) &= g((\Delta q)^N_{j,k}+\tilde q_{j,k+\frac{1}{2}})-g(\tilde q_{j,k+\frac{1}{2}}). 
\end{split}
\end{equation}
\subsection{Maximum Principle} \label{section4.3}
Having derived the semi-discrete version of our scheme in \eqref{4.65}, we next show a stability property of our scheme in the following sense. Applying our scheme to a homogeneous scalar conservation law we can illustrate the maximum principle.

We consider the maximum principle of the 2-d semi-discrete scheme \eqref{4.65} applied to homogeneous scalar conservation laws,
\begin{equation}  \label{4.67}
    \frac{d}{dt} (\Delta q)_{j,k}(t) = 
    - \frac{H^x_{j+\frac{1}{2},k}(t)-H^x_{j-\frac{1}{2},k}(t)}{\Delta x}
    - \frac{H^y_{j,k+\frac{1}{2}}(t)-H^y_{j,k-\frac{1}{2}}(t)}{\Delta y},
\end{equation}
with the numerical fluxes \eqref{4.66}.  We begin by proving the maximum principle for the deviation.

\begin{theorem}[Maximum principle]
Consider the modified scalar conservation law
\begin{equation}
    (\Delta q)_t + [f(\Delta q + \tilde q)-f(\tilde q)]_x + [g(\Delta q+\tilde q)-g(\tilde q)]_y = 0. 
\end{equation}
and the forward Euler time discretization of the 2-d scheme \eqref{4.67}
\begin{equation} \label{4.69}
    (\Delta q)^{n+1}_{j,k} 
    = (\Delta q)^n_{j,k}
    -\frac{\Delta t^n}{\Delta x} \big[ H^x_{j+\frac{1}{2},k}(t)-H^x_{j-\frac{1}{2},k}(t) \big]
    -\frac{\Delta t^n}{\Delta y} \big[ H^y_{j,k+\frac{1}{2}}(t)-H^y_{j,k-\frac{1}{2}}(t)   \big]
\end{equation}
with the numerical fluxes \eqref{4.66} and the minmod limiter \eqref{4.63}.
\\
Assume the following CFL condition holds:
 \begin{equation} \label{4.70}
    \max \left( \frac{\Delta t^n}{\Delta x} \max_{\Delta q}\left|F'(\Delta q)\right|,
    \frac{\Delta t^n}{\Delta y} \max_{\Delta q}\left|G'(\Delta q) \right| \right)\leq \frac{1}{8},
\end{equation}
then the fully-discrete scheme \eqref{4.69} satisfies the maximum principle
\begin{equation}
    \max_{j,k} \{(\Delta q)^{n+1}_{j,k}\} 
    \leq \max_{j,k} \{(\Delta q)^n_{j,k}\}.    
\end{equation}
\end{theorem}

\begin{proof}
The proof is an extension of theorem 3.1 in \cite{ref9}. We refer the reader to \ref{append.D}.
\end{proof}

With the help of Corollary 5.1 and Corollary 5.2 in \cite{ref3}, we can prove the semi-discrete scheme \eqref{4.67} satisfies the maximum principle.
\\ \\
\emph{Remark.}
Consider a general 2-d Deviation method
\begin{equation}\label{eq65} 
\begin{split}
    \frac{d (\Delta q)_{j,k}}{dt} =
    & -\frac{1}{\Delta x}\bigg[\Big(
    \hat{F}(q^L_{j+\frac{1}{2},k}, q^R_{j+\frac{1}{2},k})
    - \hat{F}(\tilde q^L_{j+\frac{1}{2},k}, \tilde q^R_{j+\frac{1}{2},k})  \Big) \\
    &\qquad\qquad\qquad- \Big(\hat{F}(q^L_{j-\frac{1}{2},k}, q^R_{j-\frac{1}{2},k})
    - \hat{F}(\tilde q^L_{j-\frac{1}{2},k}, \tilde q^R_{j-\frac{1}{2},k})    \Big)\bigg] \\
    &-\frac{1}{\Delta y}\bigg[\Big(
    \hat{G}(q^L_{j,k+\frac{1}{2}}, q^R_{j,k+\frac{1}{2}})
    - \hat{G}(\tilde q^L_{j,k+\frac{1}{2}}, \tilde q^R_{j,k+\frac{1}{2}})  \Big) \\
    &\qquad\qquad\qquad- \Big(\hat{G}(q^L_{j,k-\frac{1}{2}}, q^R_{j,k-\frac{1}{2}})
    - \hat{G}(\tilde q^L_{j,k-\frac{1}{2}}, \tilde q^R_{j,k-\frac{1}{2}})   \Big)\bigg],
\end{split}
\end{equation}
where $\hat{F}$ and $\hat{G}$ are numerical flux functions and
\begin{equation} \label{eq66}
    q^{L(R)}_{j\pm\frac{1}{2},k} = (\Delta q)^{L(R)}_{j\pm\frac{1}{2},k} +\tilde q^{L(R)}_{j\pm\frac{1}{2},k}.
\end{equation}
We find the local truncation error of the finite volume (FV) method \eqref{eq65} evaluated at the stationary solution $\tilde q$:
\begin{equation}
\begin{split}\label{eq67}
    Res(\tilde q_{j,k}) :=
    & -\frac{1}{\Delta x}\Big[
    \hat{F}(\tilde q^L_{j+\frac{1}{2},k},\tilde q^R_{j+\frac{1}{2},k})
    -\hat{F}(\tilde q^L_{j-\frac{1}{2},k},\tilde q^R_{j-\frac{1}{2},k})
    \Big] \\
    & -\frac{1}{\Delta y}\Big[
    \hat{G}(\tilde q^L_{j,k+\frac{1}{2}},\tilde q^R_{j,k+\frac{1}{2}})
    -\hat{G}(\tilde q^L_{j,k-\frac{1}{2}},\tilde q^R_{j,k-\frac{1}{2}})
    \Big] 
    =\mathcal{O}((\Delta x)^a + (\Delta y)^b).
\end{split}
\end{equation}
$N:=\min\{a,b\}$ represents the order of accuracy of the method.

Then, applying the forward Euler method in time to \eqref{eq65} and adding the stationary solution both to the left and to the right hand side of \eqref{eq65}, and using \eqref{eq67} we obtain the solution $q_{j,k}$ by
\begin{equation} \label{eq68}
\begin{split}
    q^{n+1}_{j,k} 
    = q^n_{j,k}
    &- \frac{\Delta t}{\Delta x}[
    \hat{F}(q^L_{j+\frac{1}{2},k},  q^R_{j+\frac{1}{2},k})
    - \hat{F}(q^L_{j-\frac{1}{2},k},  q^R_{j-\frac{1}{2},k})
    ] \\
    &- \frac{\Delta t}{\Delta y}[
    \hat{G}(q^L_{j,k+\frac{1}{2}},  q^R_{j,k+\frac{1}{2}})
    - \hat{G}(q^L_{j,k-\frac{1}{2}},  q^R_{j,k-\frac{1}{2}})
    ] 
    -\Delta t Res(\tilde q_{j,k}).
\end{split}
\end{equation}
Set $\hat{Q}$ is the solution of the FV method, which satisfies the following scheme
\begin{equation} \label{eq69}
\begin{split}
    \hat{Q}^{n+1}_{j,k} 
    = q^n_{j,k}
    &- \frac{\Delta t}{\Delta x}[
    \hat{F}(q^L_{j+\frac{1}{2},k},  q^R_{j+\frac{1}{2},k})
    - \hat{F}(q^L_{j-\frac{1}{2},k},  q^R_{j-\frac{1}{2},k})
    ] \\
    &- \frac{\Delta t}{\Delta y}[
    \hat{G}(q^L_{j,k+\frac{1}{2}},  q^R_{j,k+\frac{1}{2}})
    - \hat{G}(q^L_{j,k-\frac{1}{2}},  q^R_{j,k-\frac{1}{2}})
    ]. 
\end{split}
\end{equation}
Then, combining \eqref{eq69} and \eqref{eq68} yields
\begin{equation} \label{eq70}
    q^{n+1}_{j,k} = \hat{Q}^{n+1}_{j,k} -\Delta tRes(\tilde q_{j,k}).
\end{equation}
Now, we assume $\hat{Q}$ satisfies the maximum principle under some CFL conditions; i.e.,
\begin{equation} \label{eq71}
    \max_{j,k} \{\hat{Q}^{n+1}_{j,k}\} 
    \leq \max_{j,k} \{\hat{Q}^n_{j,k}\}.    
\end{equation}
According to \eqref{eq71} and the relation \eqref{eq70}, we can derive that
\begin{equation} \label{eq72}
\begin{split}
    \max_{j,k}\{ q^{n+1}_{j,k}\}
    & = \max_{j,k}\{
    \hat{Q}^{n+1}_{j,k} -\Delta tRes(\tilde q_{j,k}) \} \\
    & \leq \max_{j,k}\{\hat{Q}^{n+1}_{j,k}\} + \max_{j,k}\{\Delta tRes(\tilde q_{j,k}) \} \\
    & \leq \max_{j,k}\{\hat{Q}^n_{j,k}\} + \max_{j,k}\{\Delta tRes(\tilde q_{j,k}) \} \\
    & \leq \max_{j,k}\{
    \hat{q}^n_{j,k} +\Delta tRes(\tilde q_{j,k}) \} 
    + \max_{j,k}\{\Delta tRes(\tilde q_{j,k}) \} \\
    & \leq \max_{j,k}\{
    \hat{q}^n_{j,k}  \} + \mathcal{O}\big(\Delta t((\Delta x)^a + (\Delta y)^b)\big).
\end{split}
\end{equation}
For a large number of grid points, the error term in \eqref{eq72} is small, in fact it is of order  $N + 1$. Hence, the solution $q$ in \eqref{eq72} satisfies the maximum principle up to at most one order higher than the order of the scheme.

Since the KT-type scheme has been proved in \cite{ref9} to satisfy the maximum principle, the semi-discrete scheme constructed in section \ref{sec:3.2} (essentially) does not violate the maximum principle. 

\section{Numerical experiments and validations}\label{sec:numerics}

In this section, we present the results of a number of numerical experiments in the Euler system with gravitational source term by employing our designed well-balanced fully-discrete scheme in section \ref{sec:scheme}. In all the tests, the parameter of the $(MC-\theta)$ limiter is set as $\theta=1.5$. With the help of the CFL condition, the computation of the time-step is defined by
\begin{equation}
    \Delta t = \text{CFL} \cdot \min\bigg\{\frac{\Delta x}{\max\limits_{j,k}(a^+_{j+\frac{1}{2},k},-a^-_{j+\frac{1}{2},k})},\; 
    \frac{\Delta y}{\max\limits_{j,k}(b^+_{j,k+\frac{1}{2}},-b^-_{j,k+\frac{1}{2}})} \bigg\},
\end{equation}
with CFL (number) $=0.45$ and $\epsilon$ in the definition of local wave speeds \eqref{4.9} is given by $10^{-8}$.

Consider the 2-d Euler system with the gravitational source term as follows, 
\begin{equation} 
    \left \{
    \begin{aligned}
    & q_t +f(q)_x +g(q)_y =  S(q,x), \;\;\;  x,y\in \Omega \subset \mathbb{R},\; t>0 \\
    & q(x,y,0) = q_0(x,y), \\
    \end{aligned}
    \right.
\end{equation}
where
\begin{equation*} 
    q =
    \begin{pmatrix}
        \rho \\ \rho u_1 \\ \rho u_2 \\ \; E \;
    \end{pmatrix},\quad
    f(q) =
    \begin{pmatrix}
        \rho u_1 \\ \rho u_1^2+p \\ \rho u_1 u_2 \\ \; (E+p)u_1 \;
    \end{pmatrix},\quad
    g(q) =
    \begin{pmatrix}
        \rho u_2 \\ \rho u_1 u_2 \\ \rho u_2^2 +p \\ \; (E+p)u_2 \;
    \end{pmatrix},
\end{equation*}
and
\begin{equation*}
        S(u)=
    \begin{pmatrix}
        0 \\ -\rho\phi_x \\ -\rho\phi_y \\ \; -\rho u_1 \phi_x - \rho u_2 \phi_y \;
    \end{pmatrix}.
\end{equation*}
The given function $\phi=\phi(x,y)$ is the gravitational field and the ratio of specific heats $\gamma$ is set to be 1.4 for an ideal gas. 

\subsection{Isothermal Equilibrium}\label{5.2.1}
The first numerical experiment we consider is isothermal equilibrium presented in \cite{ref10}.

\begin{figure}[htbp]
\centering
 \begin{tabular}{ c @{\quad} c }
    \includegraphics[width=17em]{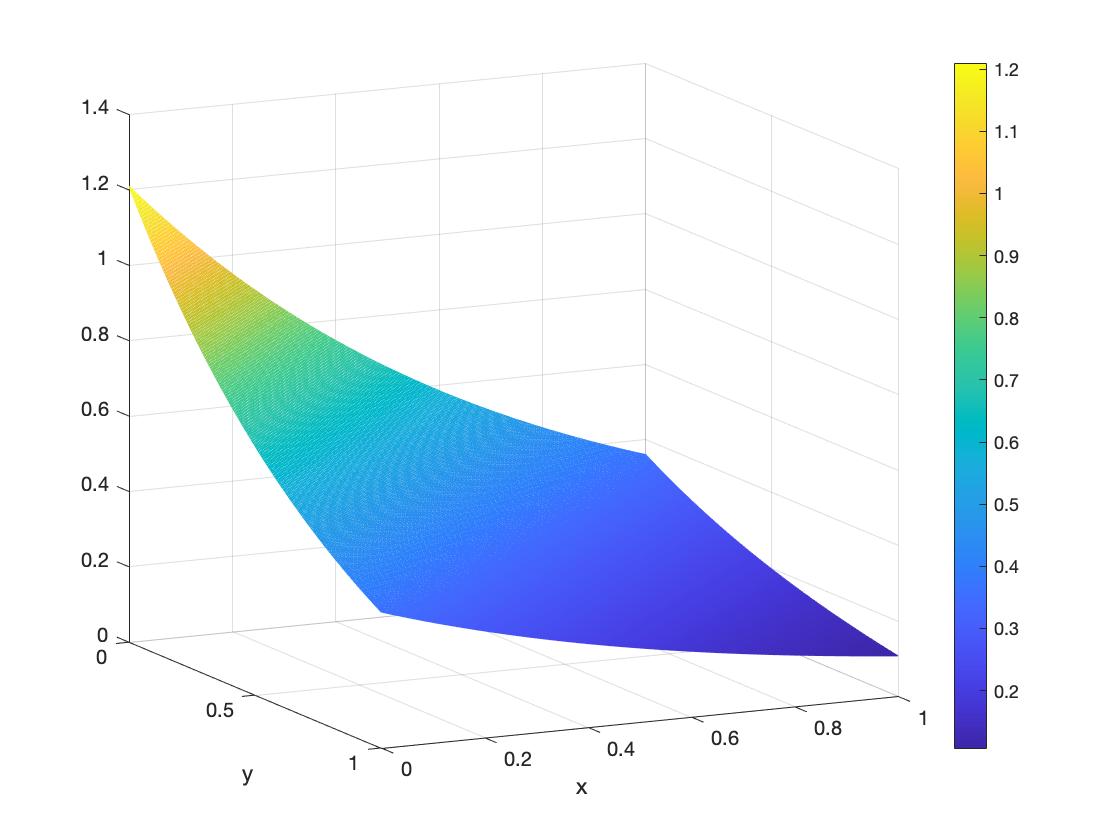} &      \includegraphics[width=17em]{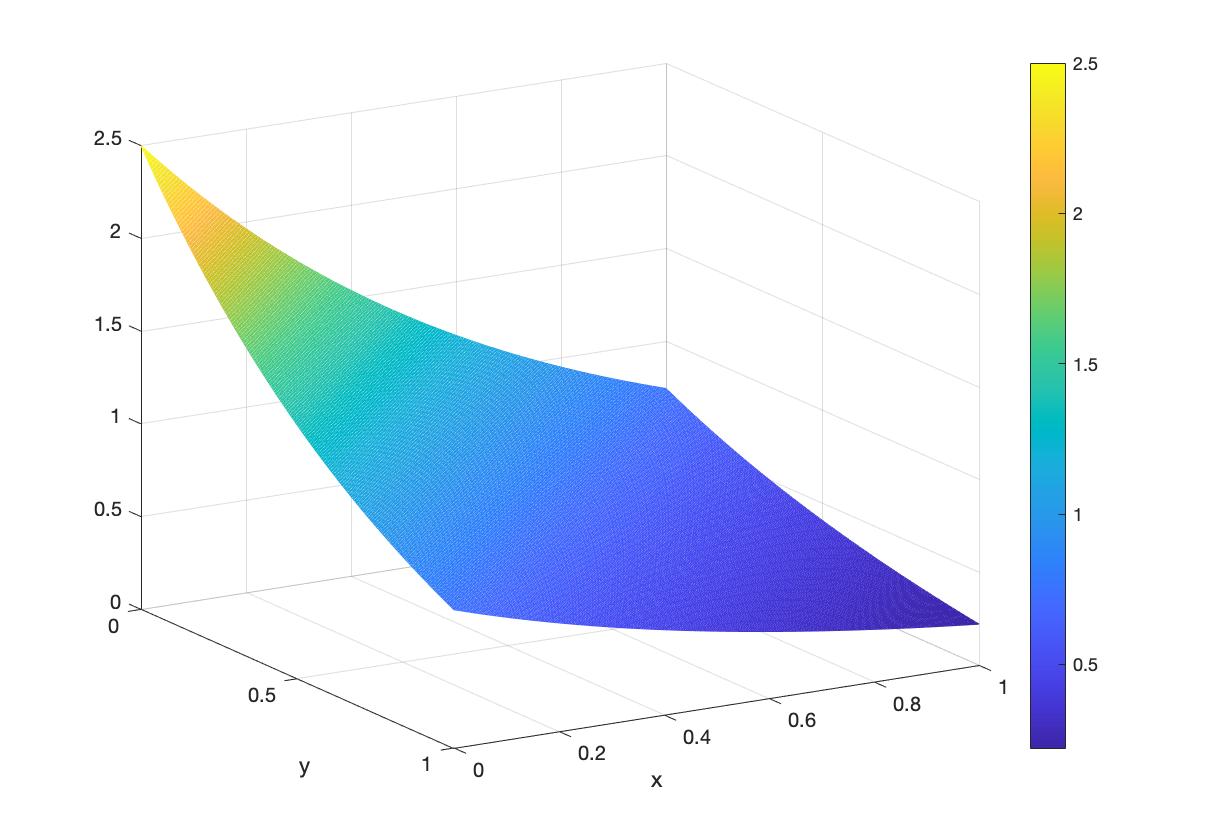} \\
    \small New scheme: Density &
      \small New scheme: Energy \\
    \includegraphics[width=17em]{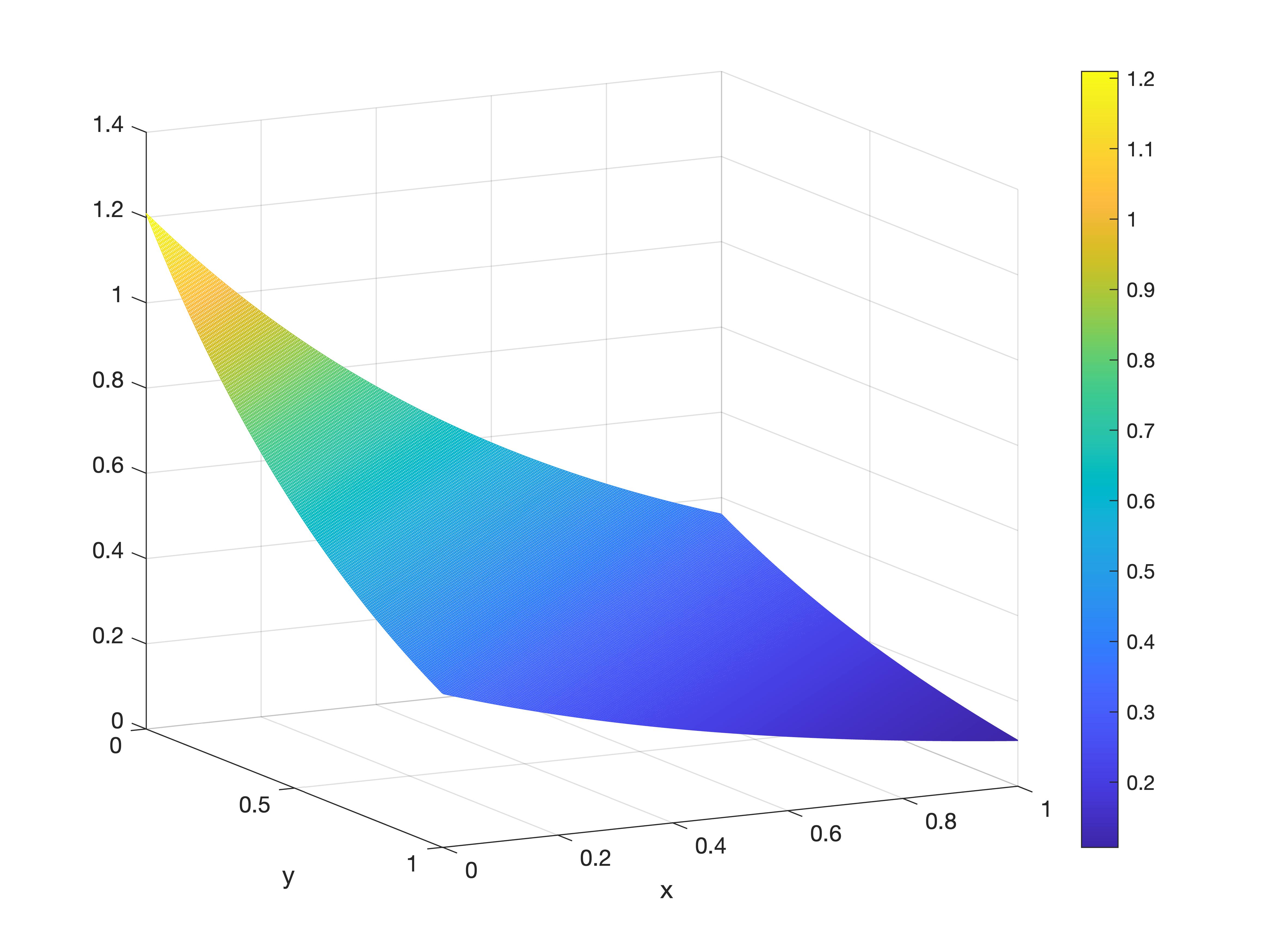} &      \includegraphics[width=17em]{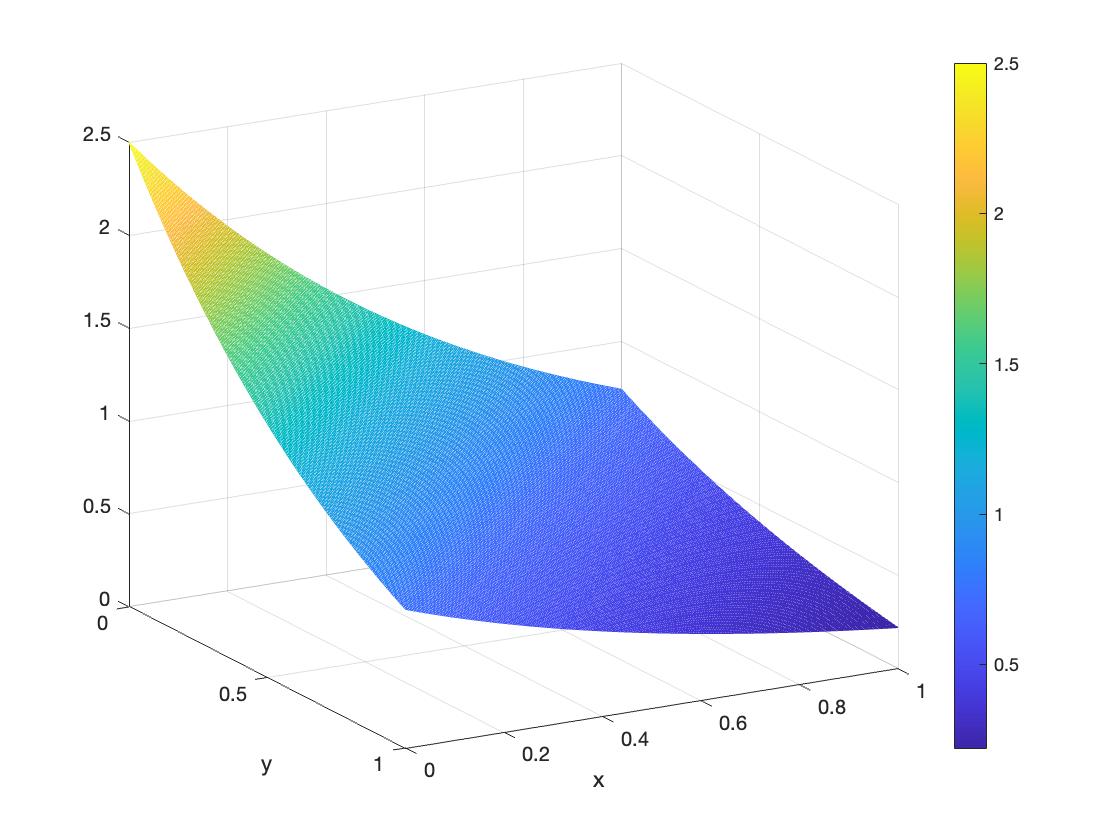} \\
    \small Exact solution:  Density &
      \small Exact solution: Energy \\
\end{tabular}
    \caption{Results of 2-d isothermal equilibrium: the top two figures are the results from our new scheme; the two below are the exact solutions.  }\label{fig.5.2-1}
\end{figure}

The 2-d isothermal equilibrium state is given by
\begin{equation} \label{5.11}
\begin{split}
    \rho(x,y) &= \rho_0 \exp(-\frac{\rho_0}{p_0}(\phi_x x + \phi_y y)),\\
    u_1(x,y) &= 0,\\
    u_2(x,y) &= 0,\\
    p(x,y) &= p_0\exp(-\frac{\rho_0}{p_0}(\phi_x x+\phi_y y)).
\end{split}
\end{equation}
Here, we set $\rho_0 = 1.21$ and $p_0 = 1$. The linear gravitational potential is given by $\phi_x = 1$ and $\phi_y = 1$. The chosen stationary solution $\tilde q$ is the isothermal equilibrium state and the outflow boundary condition is considered in this experiment. The solution is computed on the $200\times 200$ grid points in the square $[0,1]^2$ until the final time $t=0.25$ and presented in the top of the figure \ref{fig.5.2-1}. Compared to the exact solution shown in the lower part of the figure \ref{fig.5.2-1}, both of the graphics of the density and the energy are as same as the exact solutions.  

\subsection{Isothermal Equilibrium with Perturbation} \label{subsection5.2.2}
\begin{figure}[htbp]
\centering
 \begin{tabular}{ c @{\quad} c }
    \includegraphics[width=17em]{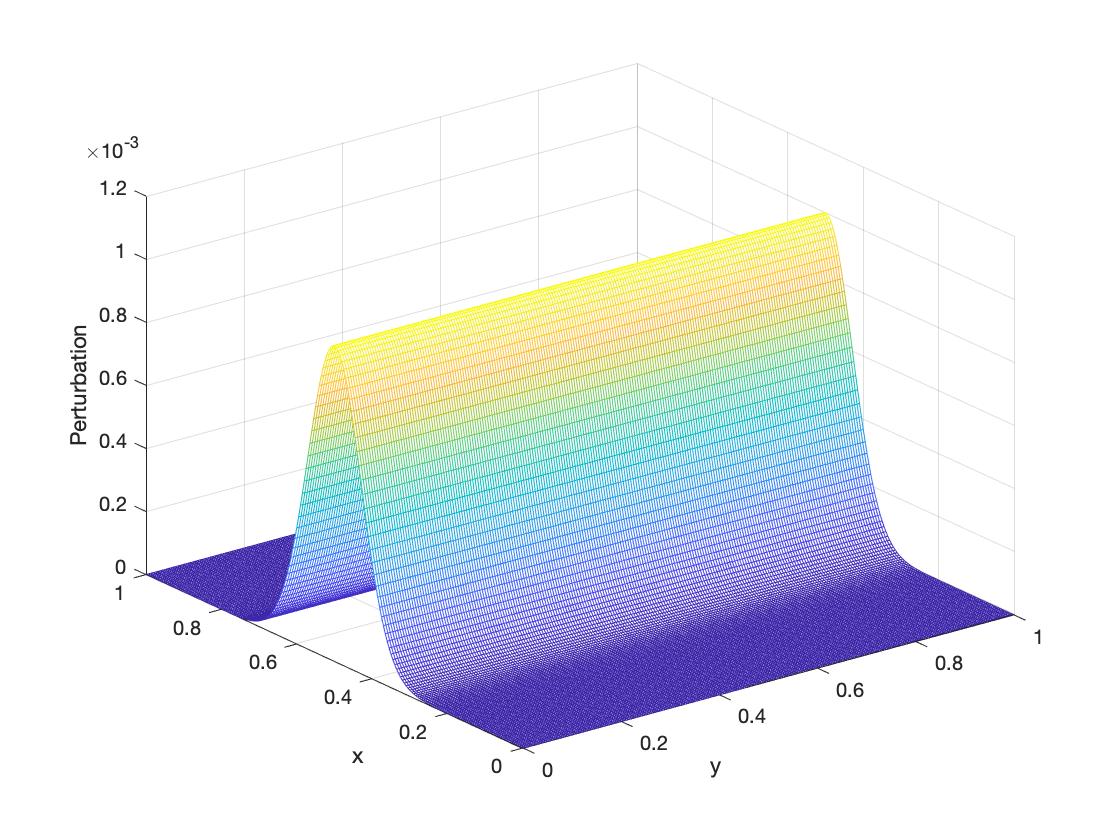} &
      \includegraphics[width=17em]{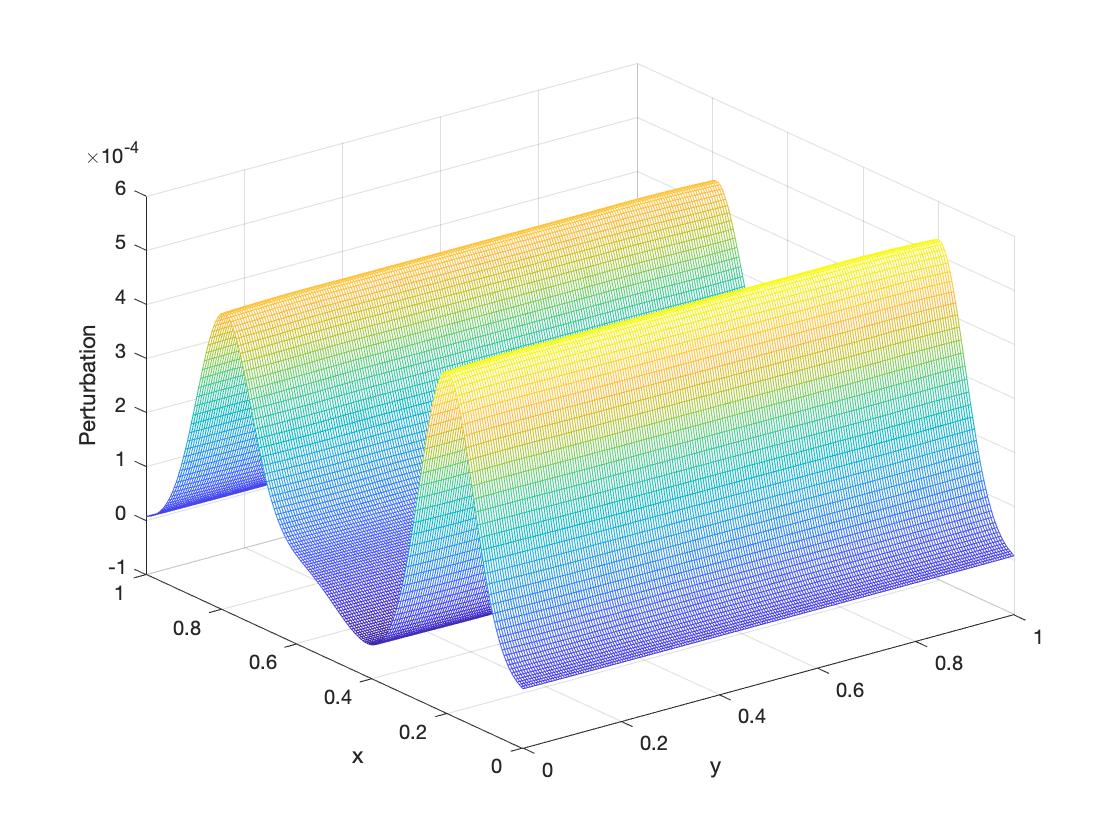} \\
    \small Initially along x &
      \small at t=0.25 along x \\
       \includegraphics[width=17em]{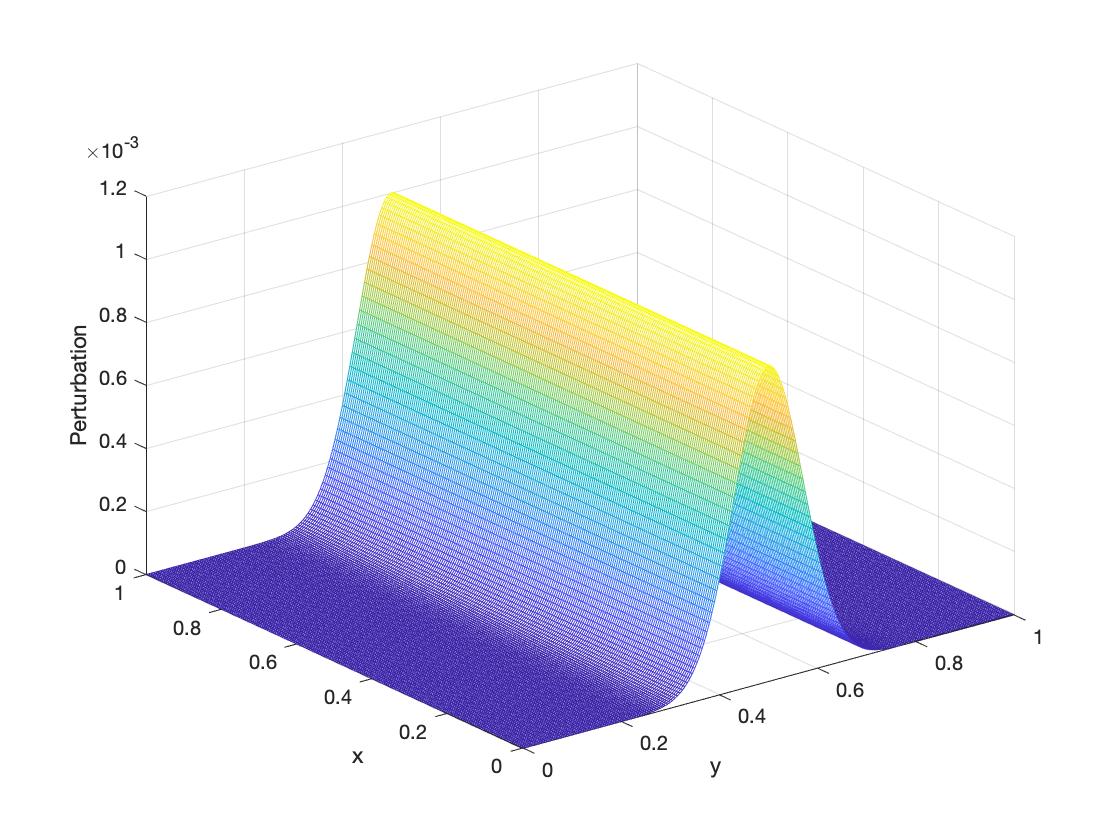} &
      \includegraphics[width=17em]{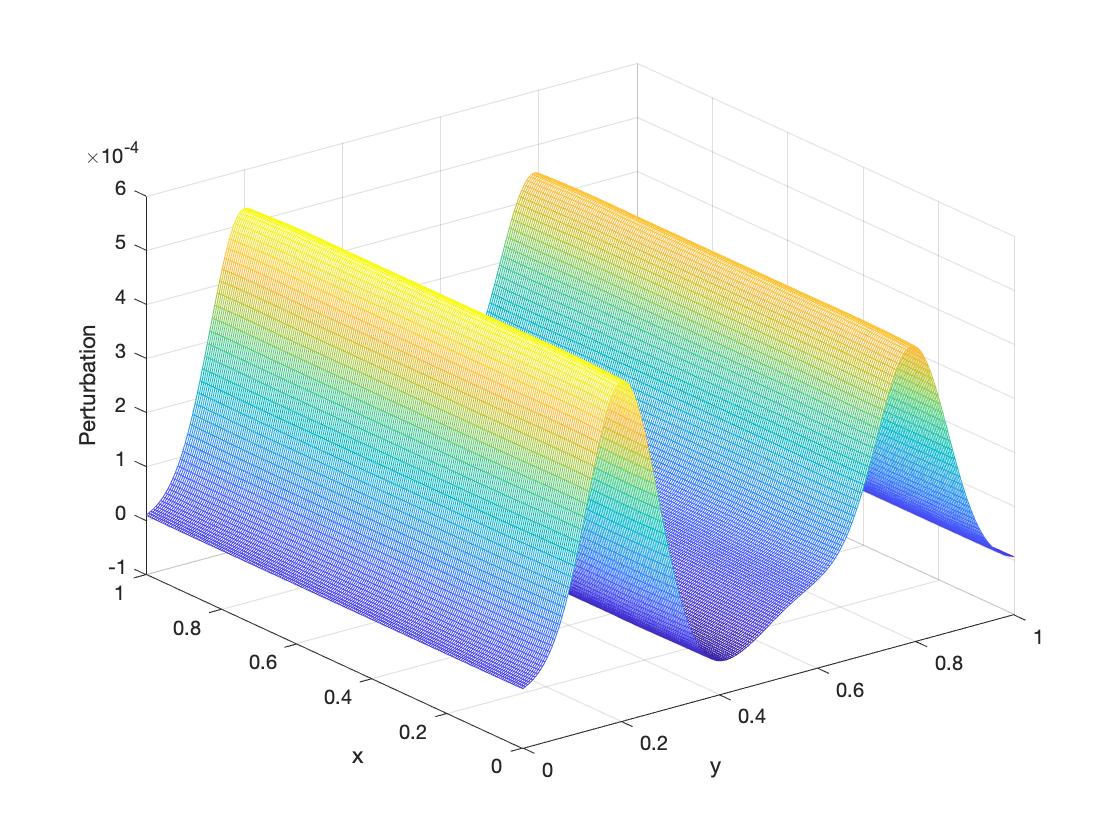} \\
    \small Initially along y &
      \small at t=0.25 along y \\
     \end{tabular}
     \caption{Results of 2-d unidirectional equilibrium perturbation.}\label{fig.5.2-2.1}
\end{figure}

\begin{figure}[htbp] 
\centering
 \begin{tabular}{ c @{\quad} c }
    \includegraphics[width=17em]{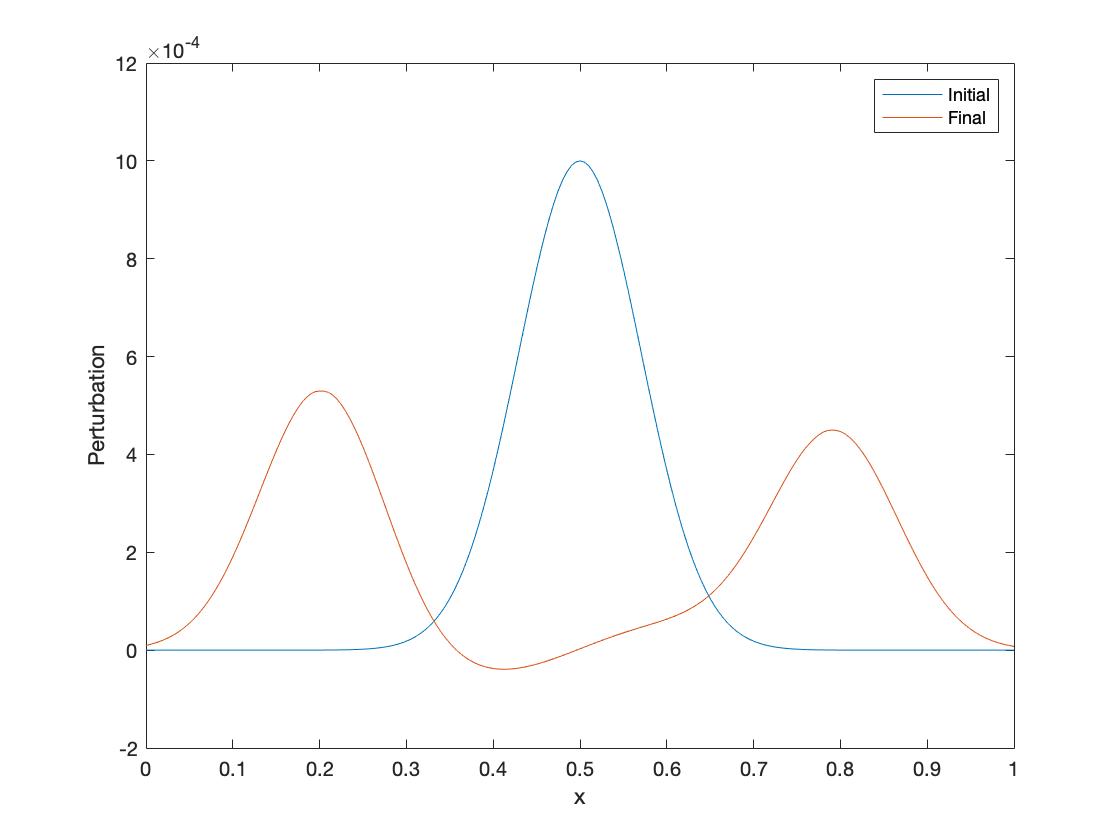} &
      \includegraphics[width=17em]{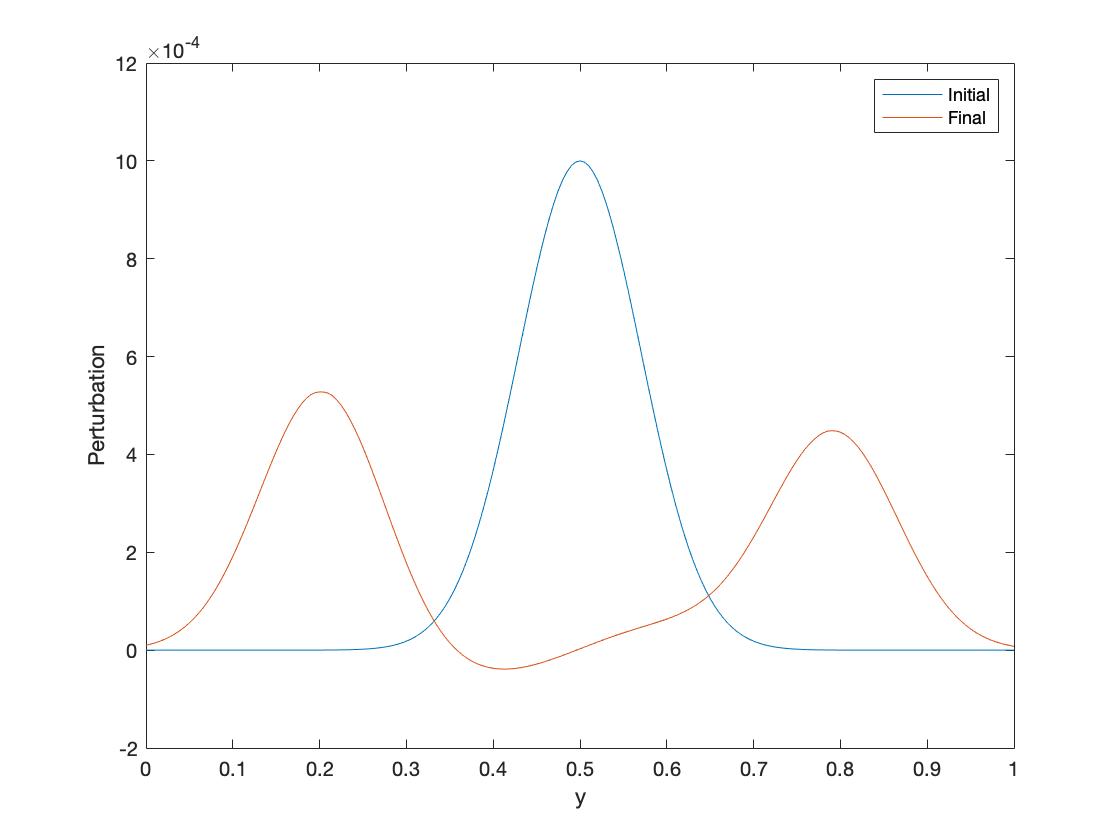} \\
    \small along x-axis &
      \small along y-axis.
 \end{tabular}
 \caption{Results of 2-d unidirectional equilibrium perturbation: Initial perturbation at t=0 and the later perturbation at the final time t=0.25.}\label{fig.5.2-2.2}
\end{figure}

This example is an extension of the previous experiment by adding a small perturbation along the $x$- or the $y$-axis to the initial pressure. For the case along $x$-axis, the initial state is set as
\begin{equation}
\begin{split}
    \rho(x,y) &= \exp(-x),\\
    u_1(x,y) &= 0,\\
    u_2(x,y) &= 0,\\
    p(x,y) &= \exp(-x) + \eta \exp(-100(x-0.5)^2).
\end{split}
\end{equation}
For the case along $y$ axis, the initial data is defined similarly. Figure \ref{fig.5.2-2.1} shows the perturbation along the $x$- and $y$- axes, respectively on $200\times200$ grid points at the initial time $t=0$ and the final time $t=0.25$. The comparison of the cross sections of the perturbation between the initial time and the final time is presented in figure \ref{fig.5.2-2.2}.

To test the order of convergence, we use the solution on a finer grid consisting of  $640\times640$ as the reference solution, and we set the CFL number to be 0.485. The results are reported in table \ref{tab.5.1}.  


\begin{table}[th]
    \centering
    \begin{tabular}[c]{|c|c|c|c|c|c|c|}
    \hline
    N & $\rho$ $L_1$-error & rate & $p$ $L_1$-error &  rate & $E$ $L_1$-error &  rate \\
    \hline
    $40^2$& 1.87E-05& - &1.67E-05 & - & 4.19E-05 & - \\
    \hline
    $80^2$& 8.93E-06 & 1.06 & 1.05E-05 & 0.67& 2.62E-05 & 0.67\\
    \hline
    $160^2$& 3.42E-06 & 1.38 & 3.89E-06 & 1.43& 9.72E-06 & 1.43\\
    \hline
    $320^2$& 1.06E-06 & 1.69 & 9.64E-07 & 2.01 & 2.41E-06 & 2.01 \\
    \hline
    \end{tabular}
    \caption{2-d isothermal equilibrium with perturbation: $L_1$-errors and convergence rates.}
    \label{tab.5.1}
\end{table}

\subsection{Moving Equilibrium} \label{subsection5.2.3}
\begin{figure}[htbp] 
\centering
 \begin{tabular}{ c @{\quad} c }
    \includegraphics[width=17em]{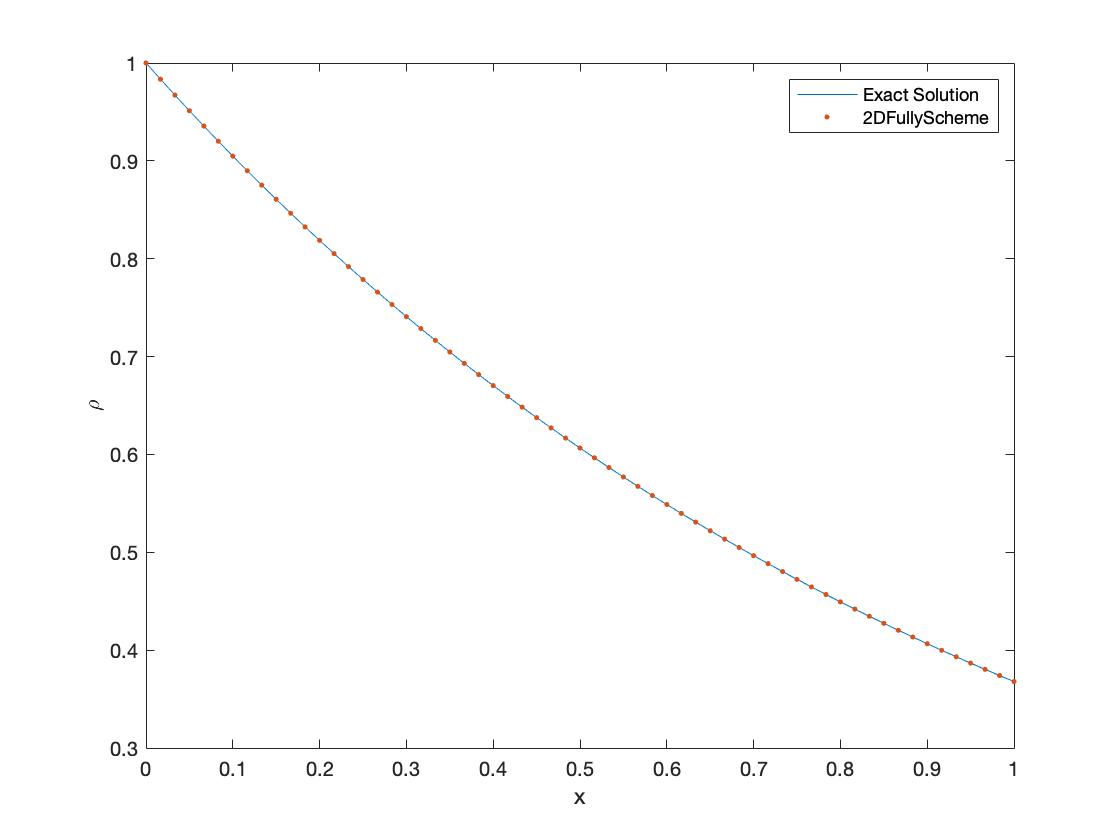} &
      \includegraphics[width=17em]{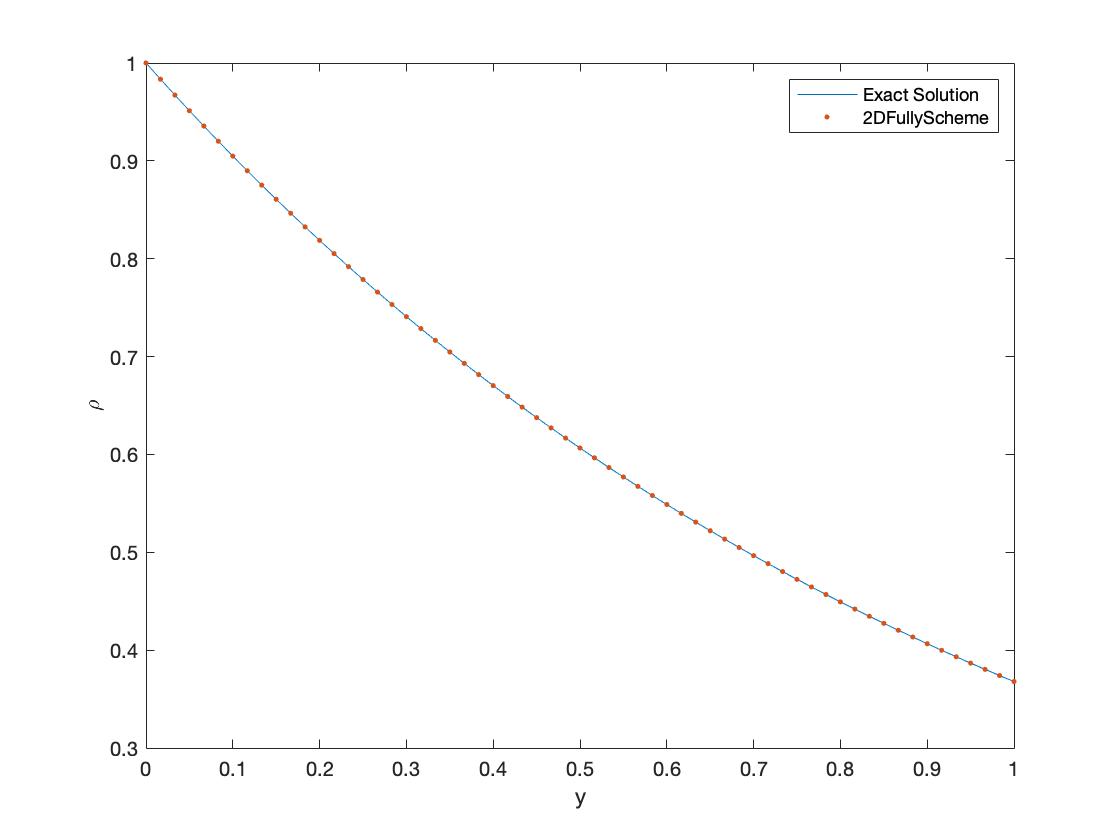} \\
    \small along x-axis &
      \small along y-axis.
      \end{tabular}
    \caption{Results of 2-d moving equilibrium at the final time t=0.25 compared to the exact solution.}\label{fig.5.2-3}
\end{figure}

Next, we test our scheme if it is capable of preserving moving equilibrium states. On the basis of the discussion in \cite{ref12}, we set the initial state of this experiment by
\begin{equation}
\begin{split}
    \rho(x,y) &= \rho_0 \exp(-\frac{\rho_0g}{p_0}(x+y)),\\
    u_1(x,y) &= \exp(x+y),\\
    u_2(x,y) &= \exp(x+y),\\
    p(x,y) &= \exp(-\frac{\rho_0g}{p_0}(x+y))^\gamma,
\end{split}
\end{equation}
with $\rho_0=1$, $p_0=1$, and $g=1$. The considered nonlinear gravitational potential is $\phi(x,y)=\exp(x+y)(-\exp(x+y)+\gamma(\exp(-\gamma(x+y))))$. We consider the experiments along the $x$- and $y$- axes separately as in the previous test case. The stationary solution used here is the equilibrium state itself and the chosen boundary condition is also the outflow boundary condition. We compute the solution on $60\times10$ grid points along the $x$-axis and on $10\times60$ grid points along the $y$-axis until the final time $t=0.25$. Figure \ref{fig.5.2-3} shows the results of cross sections along the $x$- and along the $y$- axes and we compare them to the exact solution. 

\subsection{Shock Tube Problem} \label{subsection5.2.4}
\begin{figure}[htbp]
\centering
 \begin{tabular}{ c @{\quad} c }
    \includegraphics[width=17em]{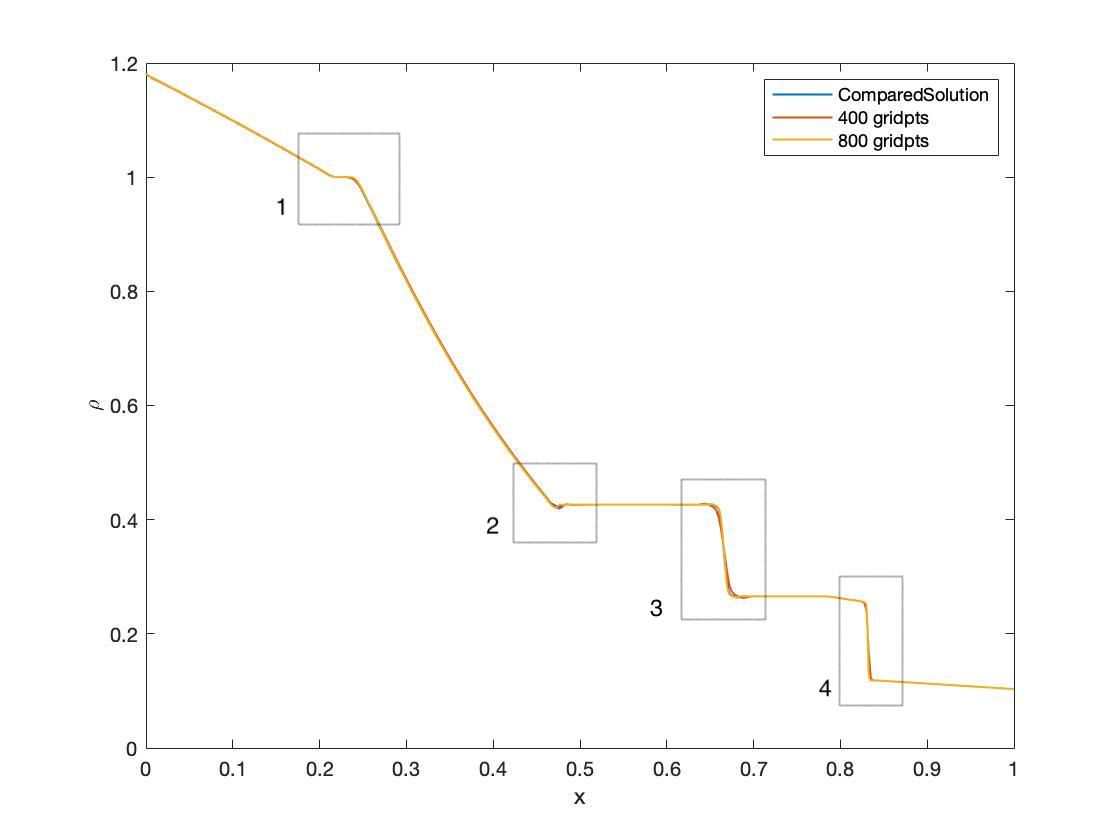} & \\
    \small Density& \\
       \includegraphics[width=17em]{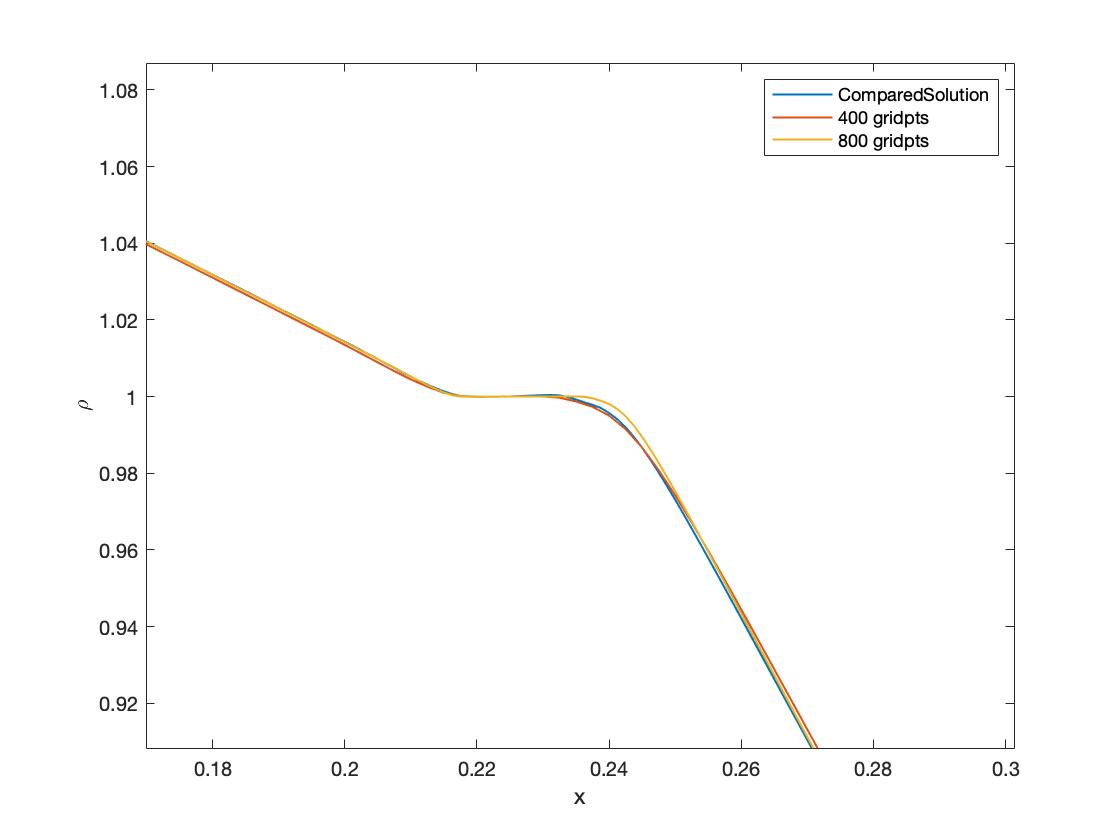}&
         \includegraphics[width=17em]{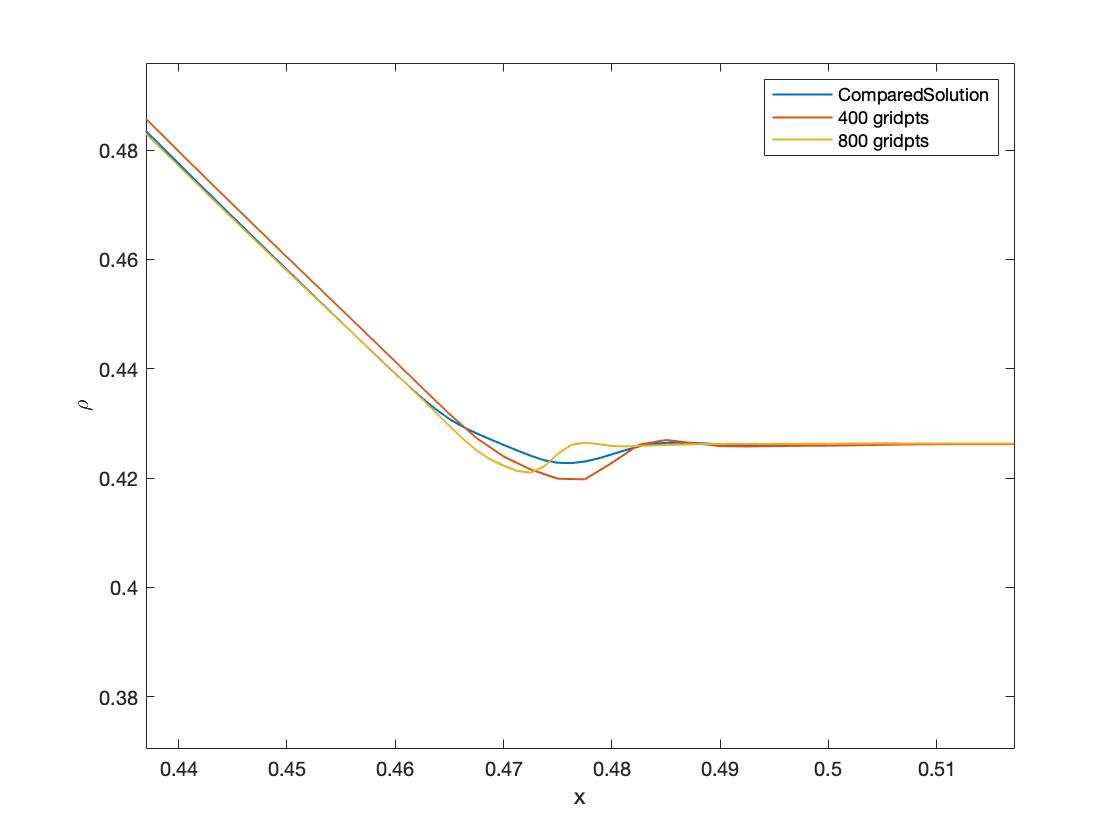}\\
    \small Zoom in on the first block&
      \small  Zoom in on the second block\\
       \includegraphics[width=17em]{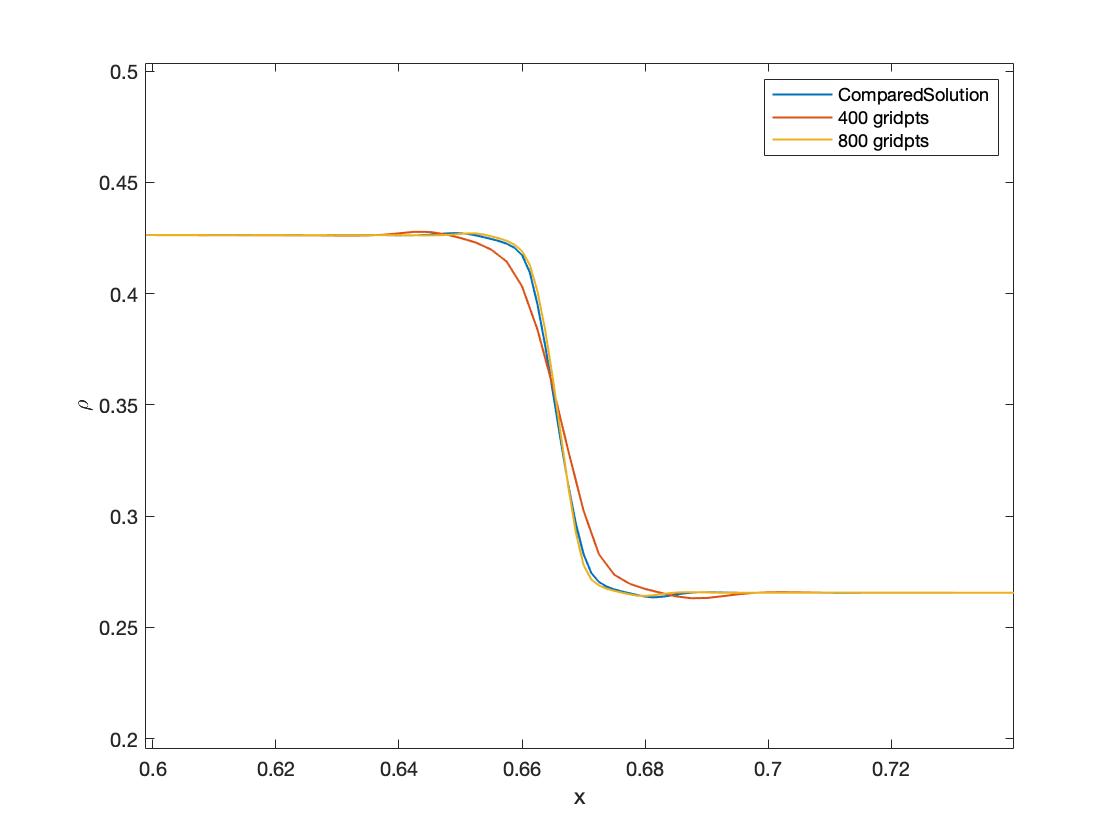}&
         \includegraphics[width=17em]{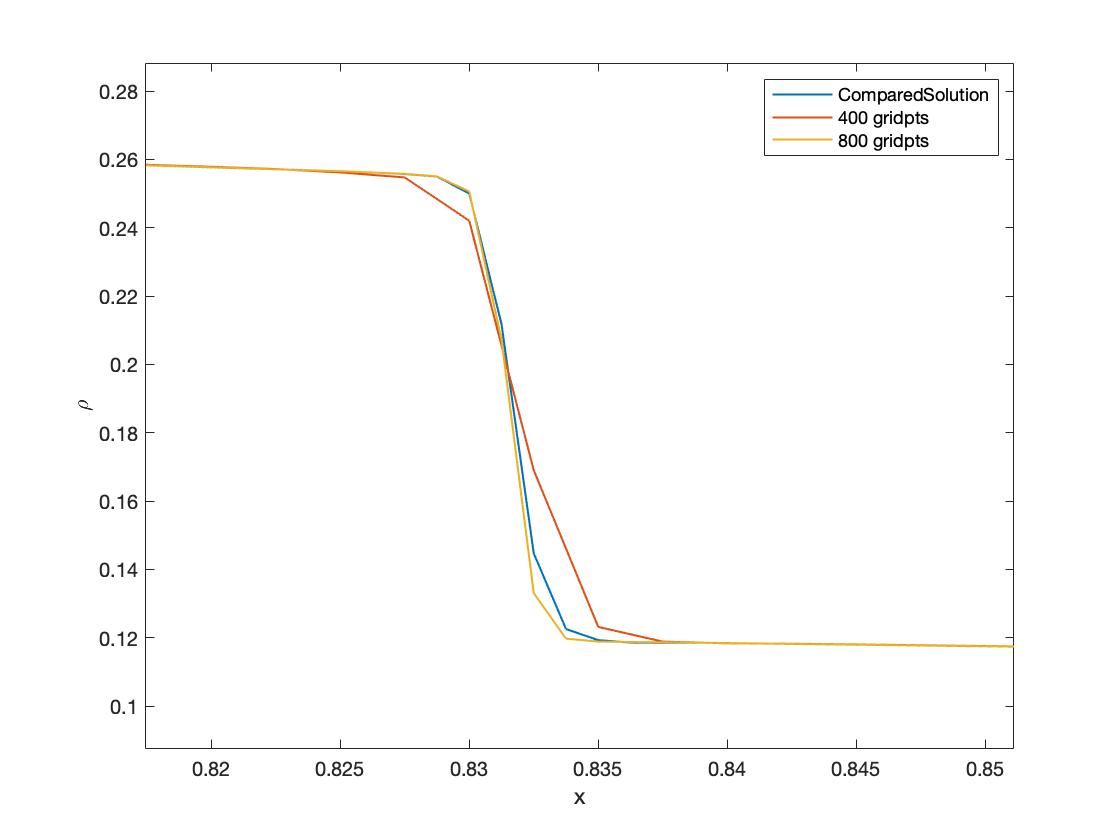}\\
    \small Zoom in on the third block&
      \small  Zoom in on the fourth block\\
\end{tabular}    
\caption{Results of 2-d shock tube problem along the $x$-axis: density and zoom in on the shocks.} \label{fig.5.2-4.1}
\end{figure}

\begin{figure}[htbp]
\centering
\begin{tabular}{ c @{\quad} c }
\includegraphics[width=17em]{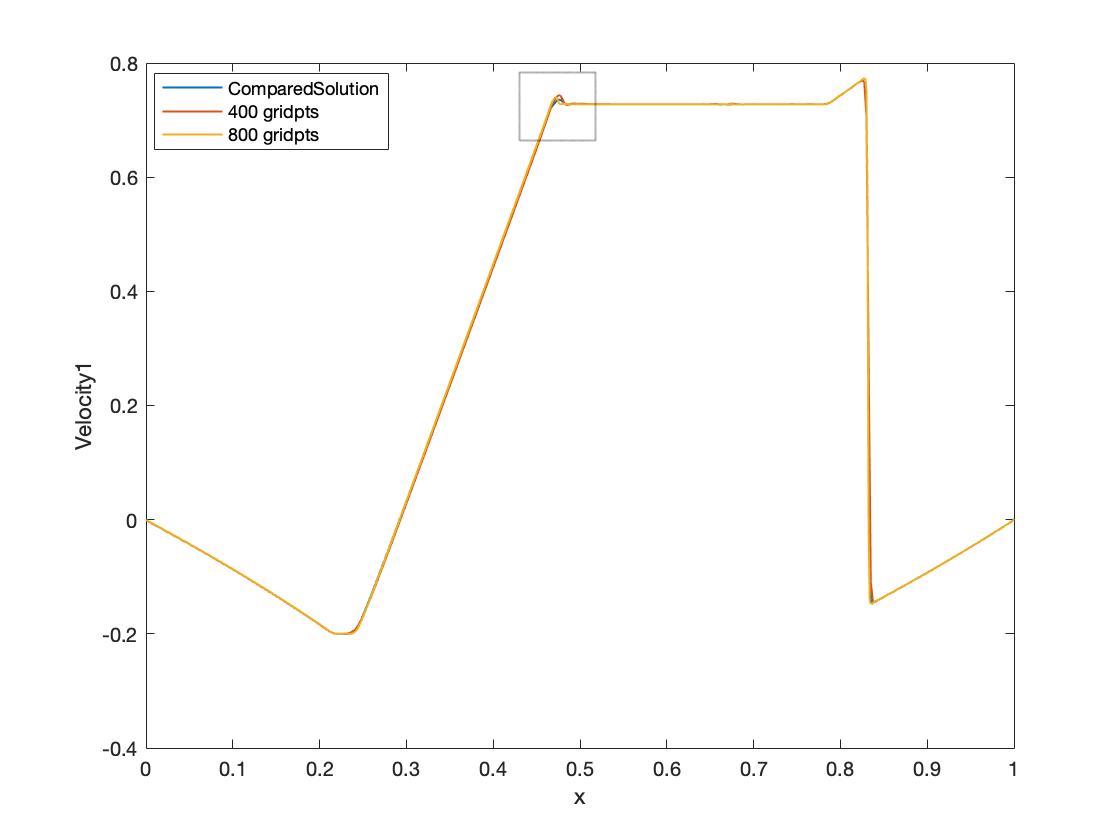} &
\includegraphics[width=17em]{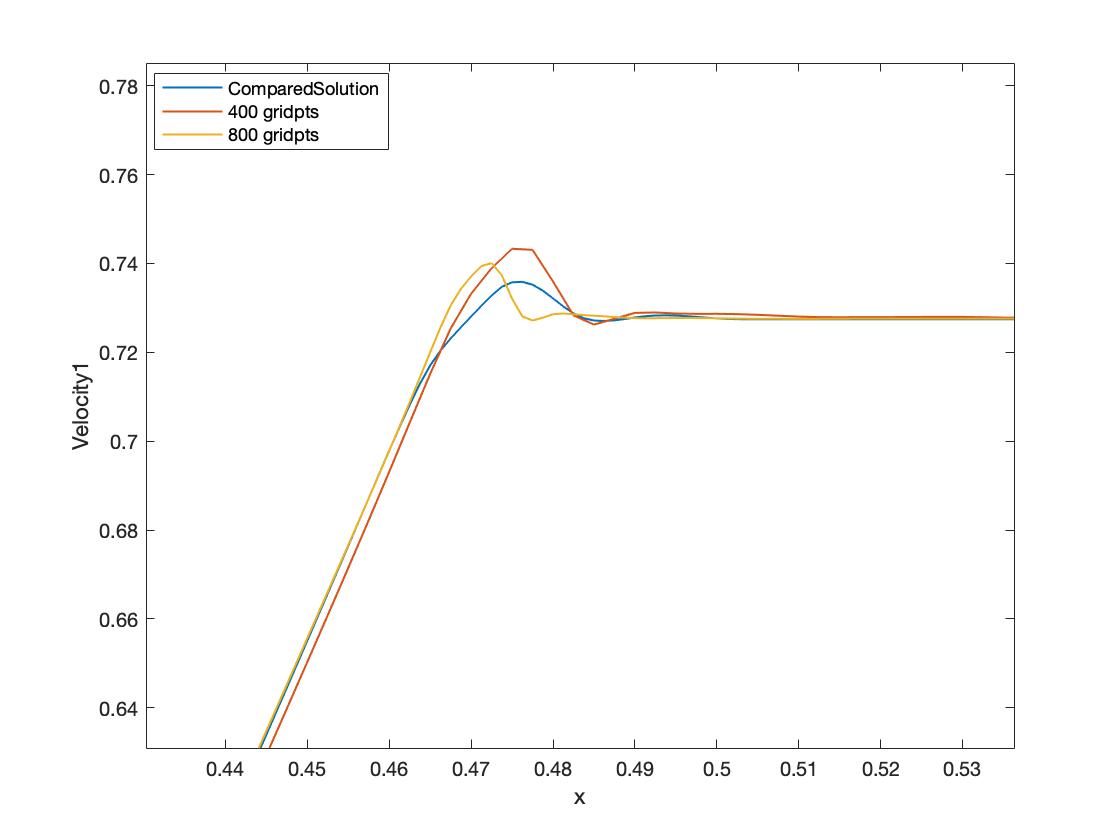} \\
\small (a) Velocity &
      \small (b) Zoom in on the block in the figure (a)\\
\includegraphics[width=17em]{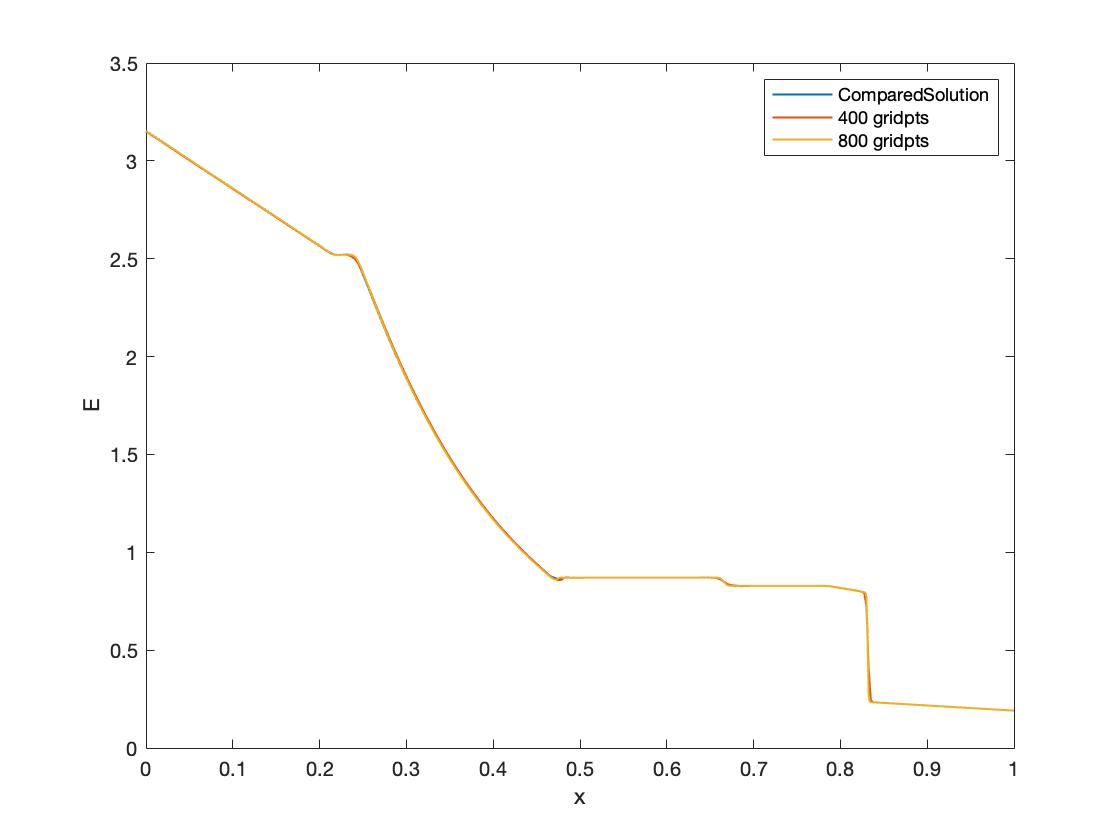} &
\includegraphics[width=17em]{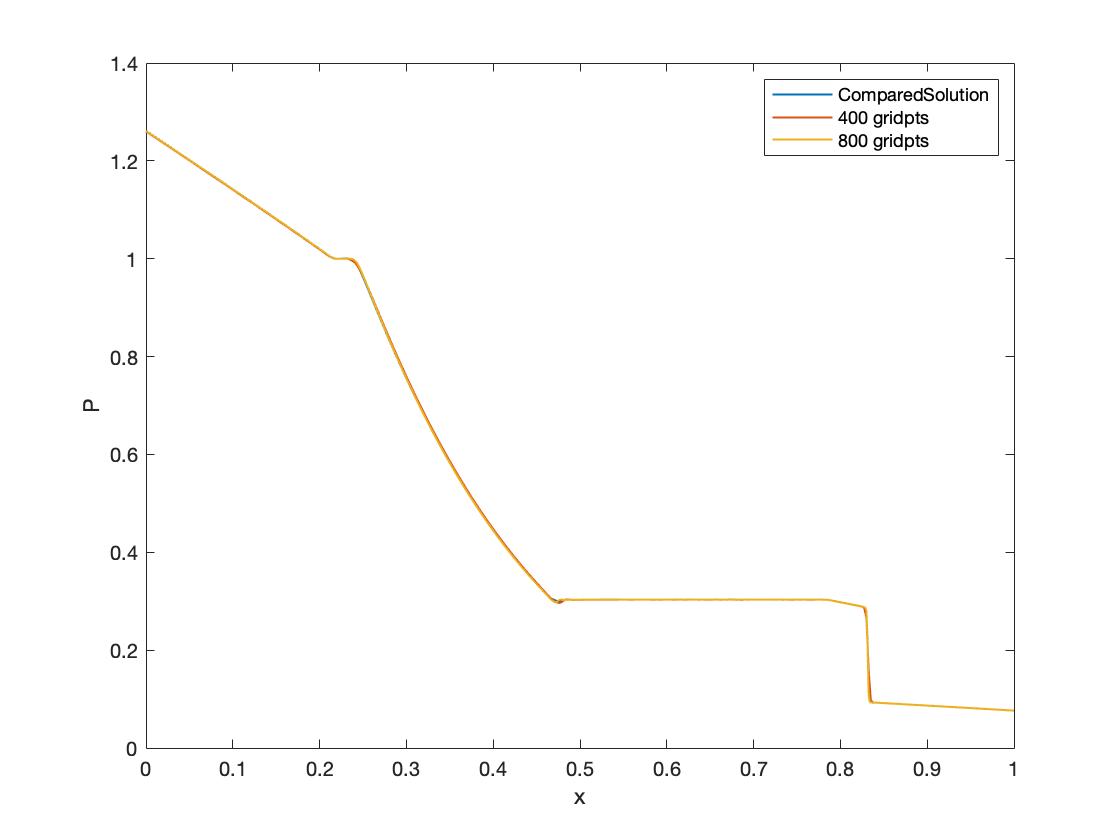} \\
\small (c) Energy &
      \small (d) Pressure.
\end{tabular}
\caption{Results of 2D shock tube problem along the $x$-axis: velocity, zoom in on the figure of velocity, energy and pressure.} \label{fig.5.2-4.2}
\end{figure}

The fourth experiment we consider is the shock tube problem along $x$-axis (see \cite{ref4}).
The initial data for this test is given by
\begin{equation}
\begin{split}
    \rho(x,y) & =
    \left \{
    \begin{aligned}
        1,\qquad\;\;\; & \text{if}\;\; x\leq0.5, \\
        0.125,\quad & \text{otherwise},
    \end{aligned}
    \right. \\
    u_1(x,y) & = 0 = u_2(x,y),\\
    p(x,y) & = 
    \left \{
    \begin{aligned}
        1,\quad\;\;\; & \text{if}\;\; x\leq0.5, \\
        0.1, \quad & \text{otherwise}.
    \end{aligned}
    \right.
\end{split}
\end{equation}
We adopt the isothermal equilibrium to be the stationary solution and choose reflecting boundary conditions for this experiment. In figures \ref{fig.5.2-4.1} and \ref{fig.5.2-4.2}, we show the solutions on $400\times 10$ and $800\times 10$ grid points obtained at the final time $t=0.2$; the reference solution on $800\times 10$ grid points is obtained following the method presented in \cite{ref4}. The performed results demonstrate our scheme can provide high resolution and the obtained solution is slightly better than the compared solution.

\begin{figure}[htbp]
\centering
          \includegraphics[width=0.45\textwidth]{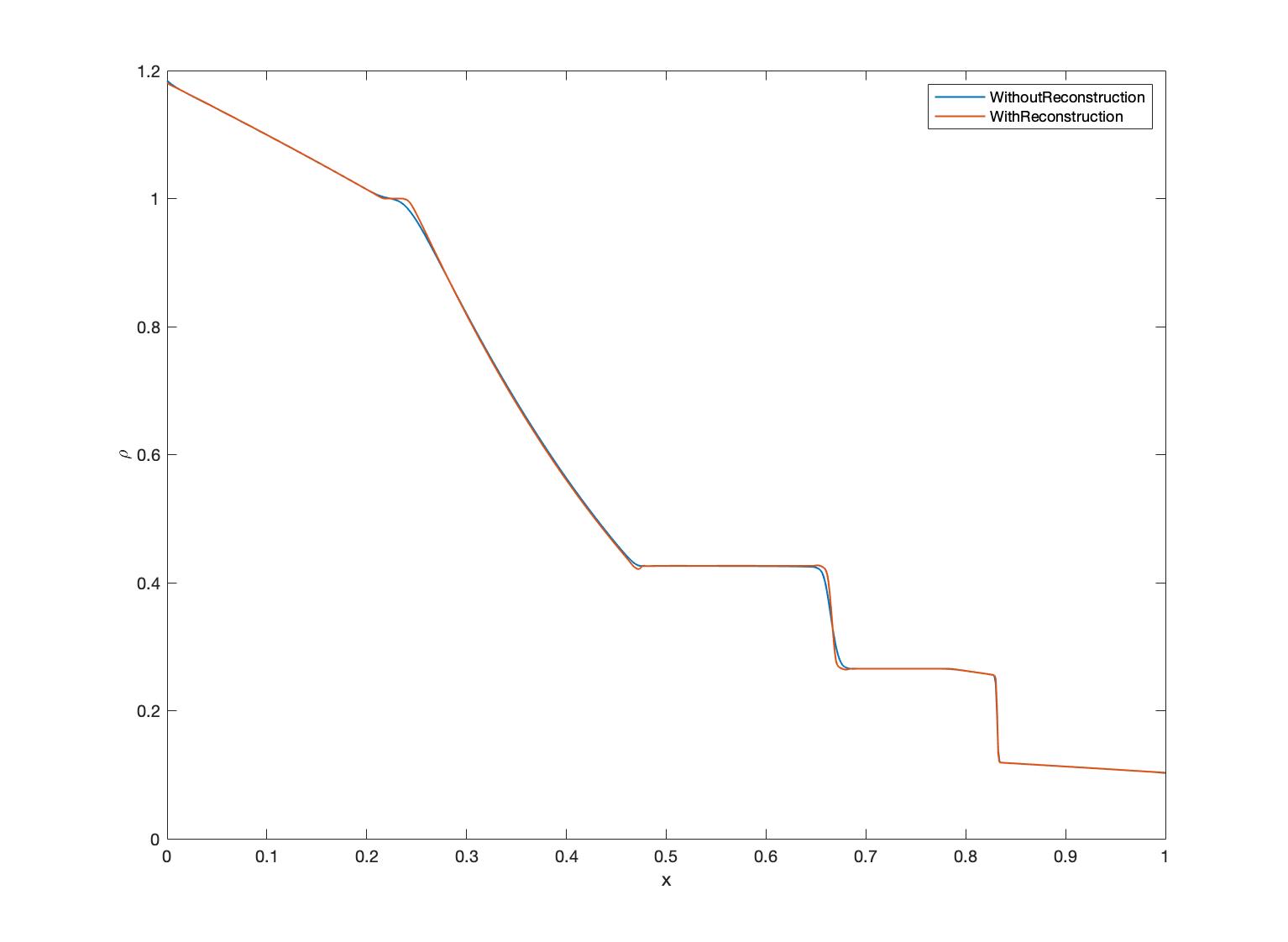}
        \caption{Results of 2D shock tube problem along the $x$-axis: profile of the density obtained with and without reconstruction.} \label{fig.5.2-4.3}
\end{figure}

In figure \ref{fig.5.2-4.3}, we compare the profile of the density on $200\times10$ grid points obtained using the reconstruction technique \eqref{4.37} at the projection step with the solution on $200\times10$ grid points obtained without reconstruction, i.e., we consider the reconstruction on non-smooth subdomains as 
\begin{equation*}
    \widetilde W^{n+1}_{D_{\alpha,\beta}} 
    = \overline{w}^{n+1}_{D_{\alpha,\beta}},
    \quad
    (x,y) \in D_{\alpha,\beta},
\end{equation*}
where $D_{\alpha,\beta}=\{D_{j+\frac{1}{2},k}, D_{j+\frac{1}{2},k+\frac{1}{2}}, D_{j,k+\frac{1}{2}}\}$. According to this figure, our reconstructed polynomial in projection step help to gain the higher resolution at the shocks.

\section{Conclusion}\label{sec:conc}
We constructed a new fully-discrete, well-balanced, central scheme by blending the ideas of the KT scheme and the Deviation method and this preserves the advantages of the mentioned two approaches. The simulations of the Euler equations with gravitational source term are closed to the reference solution and achieve the second order accuracy. For our semi-discrete scheme in section 3, it is essentially non-oscillatory from the proof of the maximum principle. Our work demonstrates that the blend of the KT scheme and the Deviation method makes a success on solving the Euler equations with gravity. In our future work, we plan to apply our scheme to other systems to examine its further practicality. 

\section*{Acknowledgement}

We thank Praveen Chandrashekar (Tata Insitute, Centre for Applicable Mathematics, Bangalore, India) who contributed to section 3 with helpful discussions.

\newcommand{\etalchar}[1]{$^{#1}$}

\appendix

\section{Approximation of source term} \label{append.A}
Based on the two facts used to approximate \eqref{4.15}, we estimate the source term on the side subdomain $D_{j+\frac{1}{2},k}$ at time $t^{n+1}$ by
\begin{equation} \label{4.30}
\begin{split}
    & \frac{1}{|D_{j+\frac{1}{2},k}|}  \int^{t^{n+1}}_{t^n} \iint_{D_{j+\frac{1}{2},k}} S(\Delta  q(x,y,t))dxdydt \\
    &\quad := \frac{\Delta t}{|D_{j+\frac{1}{2},k}|}\iint_{D_{j+\frac{1}{2},k}} S (\Delta  q(x,y,t^{n+\frac{1}{2}}))dxdy \\
    & \;\quad = \frac{\Delta t}{|D_{j+\frac{1}{2},k}|} \Bigg[
    \iint_{C^E_{j,k}} S(\Delta q(x,y,t^{n+\frac{1}{2}}))dxdy
    \\&\hspace{6cm}+ \iint_{C^W_{j+1,k}} S(\Delta q(x,y,t^{n+\frac{1}{2}}))dxdy \Bigg] \\
    & \;\quad = \frac{\Delta t}{|D_{j+\frac{1}{2},k}|} \left[
    |C^E_{j,k}| S(\Delta q(z^E_{j,k},t^{n+\frac{1}{2}}))
    + |C^W_{j+1,k}| S(\Delta q(z^W_{j+1,k},t^{n+\frac{1}{2}})) \right] \\
    &\quad\; =: \Delta t S_{D_{j+\frac{1}{2},k}}(\Delta q).
\end{split}
\end{equation}
By using the same way, the approximation on the side subdomain $D_{j,k+\frac{1}{2}}$ and the corner subdomain $D_{j+\frac{1}{2},k+\frac{1}{2}}$ can be calculated as follows,
\begin{equation} \label{4.31}
\begin{split}
    \frac{1}{|D_{j,k+\frac{1}{2}}|} &  \int^{t^{n+1}}_{t^n} \iint_{D_{j,k+\frac{1}{2}}} S(\Delta  q(x,y,t))dxdydt \\
    &\quad := \frac{\Delta t}{|D_{j,k+\frac{1}{2}}|}\iint_{D_{j,k+\frac{1}{2}}} S (\Delta  q(x,y,t^{n+\frac{1}{2}}))dxdy \\
    & \;\quad = \frac{\Delta t}{|D_{j,k+\frac{1}{2}}|} \Bigg[
    \iint_{C^N_{j,k}} S(\Delta q(x,y,t^{n+\frac{1}{2}}))dxdy
    \\&\hspace{5cm}+ \iint_{C^S_{j,k+1}} S(\Delta q(x,y,t^{n+\frac{1}{2}}))dxdy \Bigg ] \\
    &\quad := \frac{\Delta t}{|D_{j,k+\frac{1}{2}}|} \bigg[
    |C^N_{j,k}| S(\Delta q(z^N_{j,k},t^{n+\frac{1}{2}}))
    + |C^S_{j,k+1}| S(\Delta q(z^S_{j,k+1},t^{n+\frac{1}{2}})) \bigg] \\
    &\;\quad =: \Delta t S_{D_{j,k+\frac{1}{2}}}(\Delta q),
\end{split}
\end{equation}
and
\begin{equation} \label{4.32}
\begin{split}
     &\frac{1}{|D_{j+\frac{1}{2},k+\frac{1}{2}}|}  \int^{t^{n+1}}_{t^n} \iint_{D_{j+\frac{1}{2},k+\frac{1}{2}}} S(\Delta  q(x,y,t))dxdydt \\
     &\quad :=  \frac{\Delta t}{|D_{j+\frac{1}{2},k+\frac{1}{2}}|} \bigg[
    |C^{NE}_{j,k}| S(\Delta q(z^{NE}_{j,k},t^{n+\frac{1}{2}}))
    + |C^{NW}_{j+1,k}| S(\Delta q(z^{NW}_{j+1,k},t^{n+\frac{1}{2}})) \\
    & \qquad\qquad + |C^{SE}_{j,k+1}| S(\Delta q(z^{SE}_{j,k+1},t^{n+\frac{1}{2}}))
    + |C^{SW}_{j+1,k+1}| S(\Delta q(z^{SW}_{j+1,k+1},t^{n+\frac{1}{2}})) \bigg] \\
    &\quad\; =: \Delta t S_{D_{j+\frac{1}{2},k+\frac{1}{2}}}(\Delta q)
\end{split}
\end{equation}
Similar to the approximation \eqref{4.25}, the value at the center of mass, $\Delta q(z^I_{j,k},t^{n+\frac{1}{2}})$, is given by the Taylor expansion, 
\begin{equation} \label{4.33}
\begin{split}
    \Delta q^{n+\frac{1}{2}}_{z^I_{j,k}} := \Delta q^n_{z^I_{j,k}} + & \frac{\Delta t}{2}(\Delta q_t)^n_{z^I_{j,k}} \\
    = \Delta q^n_{z^I_{j,k}} - & \frac{\Delta t}{2}\left[ \left[ (f(\Delta q+\tilde q)-f(\tilde q))_x \right]^n_{z^I_{j,k}} \right.\\ 
    & \quad\; +\left. \left[ (g(\Delta q+\tilde q)-g(\tilde q))_y \right]^n_{z^I_{j,k}}\right] + \frac{\Delta t}{2} S(\Delta q)^n_{z^I_{j,k}} \\
    := \Delta q^n_{z^I_{j,k}} - & \frac{\Delta t}{2}\left[ \left[ (f(\Delta q+\tilde q)-f(\tilde q))_x \right]^n_{j,k} 
     + \left[ (g(\Delta q+\tilde q)-g(\tilde q))_y \right]^n_{j,k}\right] \\
     & \qquad\qquad\qquad\qquad\qquad\qquad\qquad\quad + \frac{\Delta t}{2} S(\Delta q)^n_{z^I_{j,k}}. 
\end{split}
\end{equation}

\section{Intermediate values at the evolution step in section 3.1.2} \label{append.B}

Applying \eqref{4.14} to the rectangle subdomains in figure \ref{Fig.4.5} yield the intermediate values on these subdomains at time $t^{n+1}$,  
\begin{equation} \label{4.45}
\begin{split}
    \overline{w}^{n+1}_{j+\frac{1}{2},k}& = \frac{1}{|D_{j+\frac{1}{2},k}|} 
    \left [ \iint_{D_{j+\frac{1}{2},k}} Q(x,y,t^n) \,dxdy \right. \\
    & - \int^{t^{n+1}}_{t^n} \int^{y_{k+\frac{1}{2}}+B^-_{j+\frac{1}{2},k+\frac{1}{2}}\Delta t}_{y_{k-\frac{1}{2}}+B^+_{j+\frac{1}{2},k-\frac{1}{2}}\Delta t} \Big[f(\Delta q +\tilde q) -f(\tilde q)\Big]^{x_{j+\frac{1}{2}}+a^+_{j+\frac{1}{2},k}\Delta t}_{x=x_{j+\frac{1}{2}}+a^-_{j+\frac{1}{2},k}\Delta t}  \,dydt \\
    & - \int^{t^{n+1}}_{t^n} \int^{x_{j+\frac{1}{2}}+a^+_{j+\frac{1}{2},k}\Delta t}_{x_{j+\frac{1}{2}}+a^-_{j+\frac{1}{2},k}\Delta t} \Big[g(\Delta q +\tilde q) -g(\tilde q)\Big]^{y_{j+\frac{1}{2}}+B^-_{j+\frac{1}{2},k+\frac{1}{2}}\Delta t}_{y=y_{k-\frac{1}{2}}+B^+_{j+\frac{1}{2},k+\frac{1}{2}}\Delta t} \,dxdt \\
    &\qquad\;\; \left. + \int^{t^{n+1}}_{t^n}\iint_{D_{j+\frac{1}{2},k}} S(\Delta q) \,dxdydt \right ],    
\end{split}
\end{equation}
\begin{equation}  \label{4.46}
\begin{split}
    \overline{w}^{n+1}_{j,k+\frac{1}{2}} =& \frac{1}{|D_{j,k+\frac{1}{2}}|} 
    \left [ \iint_{D_{j,k+\frac{1}{2}}} Q(x,y,t^n) \,dxdy \right. \\
    & - \int^{t^{n+1}}_{t^n} \int^{y_{k+\frac{1}{2}}+b^+_{j,k+\frac{1}{2}}\Delta t}_{y_{k+\frac{1}{2}}+b^-_{j,k+\frac{1}{2}}\Delta t} \Big[f(\Delta q +\tilde q) -f(\tilde q)\Big]^{x_{j+\frac{1}{2}}+A^-_{j+\frac{1}{2},k+\frac{1}{2}}\Delta t}_{x=x_{j-\frac{1}{2}}+A^+_{j-\frac{1}{2},k+\frac{1}{2}}\Delta t}  \,dydt \\
    & - \int^{t^{n+1}}_{t^n} \int^{x_{j+\frac{1}{2}}+A^-_{j+\frac{1}{2},k+\frac{1}{2}}\Delta t}_{x_{j-\frac{1}{2}}+A^+_{j-\frac{1}{2},k+\frac{1}{2}}\Delta t} \Big[g(\Delta q +\tilde q) -g(\tilde q)\Big]^{y_{k+\frac{1}{2}}+b^+_{j,k+\frac{1}{2}}\Delta t}_{y=y_{k+\frac{1}{2}}+b^-_{j,k+\frac{1}{2}}\Delta t} \,dxdt \\
    &\qquad \left. + \int^{t^{n+1}}_{t^n}\iint_{D_{j,k+\frac{1}{2}}} S(\Delta q) \,dxdydt \right ],    
\end{split}
\end{equation}
\begin{equation} \label{4.47}
\begin{split}
    \overline{w}^{n+1}_{j+\frac{1}{2},k+\frac{1}{2}} = \frac{1}{|D_{j+\frac{1}{2},k+\frac{1}{2}}|} 
    &\left [ \iint_{D_{j+\frac{1}{2},k+\frac{1}{2}}} Q(x,y,t^n) \,dxdy \right. \\
    &  - \int^{t^{n+1}}_{t^n} \iint_{D_{j+\frac{1}{2},k+\frac{1}{2}}} \Big[f(\Delta q+\tilde q)-f(\tilde q)\Big]_x \,dxdydt \\
    & - \int^{t^{n+1}}_{t^n} \iint_{D_{j+\frac{1}{2},k+\frac{1}{2}}} \Big[g(\Delta q+\tilde q)-g(\tilde q)\Big]_y \,dxdydt \\
    &\left. + \int^{t^{n+1}}_{t^n}\iint_{D_{j+\frac{1}{2},k+\frac{1}{2}}} S(\Delta q) \,dxdydt \right ],    
\end{split}
\end{equation}
and
\begin{equation} \label{4.48}
\begin{split}
   \overline{w}^{n+1}_{j,k} = \frac{1}{|D_{j,k}|} 
    &\left [ \iint_{D_{j,k}} Q(x,y,t^n) \,dxdy \right. \\
    &- \int^{t^{n+1}}_{t^n} \iint_{D_{j,k}} [f(\Delta q+\tilde q)-f(\tilde q)]_x \,dxdydt \\
    & - \int^{t^{n+1}}_{t^n} \iint_{D_{j,k}} [g(\Delta q+\tilde q)-g(\tilde q)]_y \,dxdydt \\
    & \left. + \int^{t^{n+1}}_{t^n}\iint_{D_{j,k}} S(\Delta q) \,dxdydt \right ].    
\end{split}
\end{equation}

\section{Derivation of the semi-discrete scheme} \label{append.C}
Due to the conservation property, the last integral on the right-hand-side of \eqref{4.52} has the relation that
\begin{equation} \label{4.53}
    \iint_{D_{j,k}} \overline{w}^{n+1}_{j,k}(x,y) dxdy = |D_{j,k}|\overline{w}^{n+1}_{j,k}.
\end{equation}
For $S = {(j\pm\frac{1}{2},k), (j,k\pm\frac{1}{2}), (j\pm\frac{1}{2},k\pm\frac{1}{2})}$, according to the approximation of $\widetilde W^{n+1}_S(x,y)$ in \eqref{4.50}, we have the relation between $\widetilde W^{n+1}_S$ and $\overline{w}^{n+1}_S$:
\begin{equation} \label{4.54}
     \widetilde W^{n+1}_S(x,y) = \overline{w}^{n+1}_S + O(\Delta t).
\end{equation}
Since the areas of the domain $C^{E(W)}_{j,k}$, $C^{N(S)}_{j,k}$, and $C^{NE(NW)(SE)(SW)}_{j,k}$ are evaluated by
\begin{equation} \label{4.55}
\begin{split}
    |C^{E(W)}_{j,k}| & = \Delta t\Delta y (\mp a^\mp_{j\pm\frac{1}{2},k})+O((\Delta t)^2), \\
    |C^{N(S)}_{j,k}| & = \Delta t\Delta x (\mp b^\mp_{j,k\pm\frac{1}{2}})+O((\Delta t)^2), \\
    |C^{NE(NW)(SE)(SW)}_{j,k}| & = O((\Delta t)^2),
\end{split}
\end{equation}
the relation between $\widetilde W^{n+1}_S$ and $\overline{w}^{n+1}_S$ can be written as
\begin{equation} \label{4.56}
\begin{split}
    & \iint_{C^{E(W)}_{j,k}} \widetilde W_{j\pm\frac{1}{2},k} dxdy =  
    |C^{E(W)}_{j,k}|\overline{w}^{n+1}_{j\pm\frac{1}{2},k} +O((\Delta t)^2), \\
    & \iint_{C^{N(S)}_{j,k}} \widetilde W_{j,k\pm\frac{1}{2}} dxdy = 
    |C^{N(S)}_{j,k}|\overline{w}^{n+1}_{j,k\pm\frac{1}{2}} +O((\Delta t)^2), \\
    & \iint_{C^{NE(NW)(SE)(SW)}_{j,k}} \widetilde W^{n+1}_{j\pm\frac{1}{2},k\pm\frac{1}{2}} dxdy =  O((\Delta t)^2).
\end{split}
\end{equation}
As $\Delta t \to 0$, the values at the corner vanish, because the area are proportional to $(\Delta t)^2$.\quad Thus, \eqref{4.52} reduces to
\begin{equation} \label{4.57}
\begin{split}
    \frac{d(\Delta q)_{j,k}(t)}{dt} = \lim_{\Delta t\to 0}
    \frac{1}{\Delta t} \bigg [ \frac{1}{\Delta x \Delta y} \Big (
    |C^{E(W)}_{j,k}|\overline{w}^{n+1}_{j\pm\frac{1}{2},k} 
    + |C^{N(S)}_{j,k}| \overline{w}^{n+1}_{j,k\pm\frac{1}{2}}
    \\+ |D_{j,k}|\overline{w}^{n+1}_{j,k} \Big ) 
    - \Delta q^n_{j,k} \bigg ]
\end{split}
\end{equation}
Applying the areas of $C^{E(W)}_{j,k}$ and $C^{N(S)}_{j,k}$ to the first and the second term on the right-hand-side of \eqref{4.57} respectively, we obtain that
\begin{equation} \label{4.58}
\begin{split}
    \lim_{\Delta t\to 0}\frac{|C^{E(W)}_{j,k}|}{\Delta t \Delta x \Delta y}\overline{w}^{n+1}_{j\pm\frac{1}{2},k} =
    -\frac{a^\mp_{j\pm\frac{1}{2},k}}{\Delta x} \lim_{\Delta t\to 0} \overline{w}^{n+1}_{j\pm\frac{1}{2},k}, \\
    \lim_{\Delta t\to 0}\frac{|C^{N(S)}_{j,k}|}{\Delta t \Delta x \Delta y}\overline{w}^{n+1}_{j,k\pm\frac{1}{2}} =
    -\frac{b^\mp_{j,k\pm\frac{1}{2}}}{\Delta y} \lim_{\Delta t\to 0} \overline{w}^{n+1}_{j,k\pm\frac{1}{2}}.
    \end{split}
\end{equation}
Then using the approximations of $\overline{w}^{n+1}_{j+\frac{1}{2},k}$ and $\overline{w}^{n+1}_{j,k+\frac{1}{2}}$in \eqref{4.45} and \eqref{4.46} to substitute the values in \eqref{4.58} results in
\begin{equation} \label{4.59}
\begin{split}
    \lim_{\Delta t\to 0}\frac{|C^{E(W)}_{j,k}|}{\Delta t \Delta x \Delta y}\overline{w}^{n+1}_{j\pm\frac{1}{2},k} =
    & - \frac{1}{(a^+_{j\pm\frac{1}{2},k}-a^-_{j\pm\frac{1}{2},k})\Delta x} \times \\
    &\qquad\qquad \left [
     a^-_{j\pm\frac{1}{2},k}a^+_{j\pm\frac{1}{2},k} (\Delta q)^{W(E)}_{j\pm1,k}-(a^\mp_{j\pm\frac{1}{2},k})^2 (\Delta q)^{E(W)}_{j,k}\right] \\
    & + \frac{a^\mp_{j\pm\frac{1}{2},k}}{(a^+_{j\pm\frac{1}{2},k}-a^-_{j\pm\frac{1}{2},k})\Delta x} \left[ F((\Delta q)^{W(E)}_{j\pm1,k})-F((\Delta q)^{E(W)}_{j,k}) \right], \\
    \lim_{\Delta t\to 0}\frac{|C^{N(S)}_{j,k}|}{\Delta t \Delta x \Delta y}\overline{w}^{n+1}_{j,k\pm\frac{1}{2}} =
    & - \frac{1}{(b^+_{j,k\pm\frac{1}{2}}-b^-_{j,k\pm\frac{1}{2}})\Delta y} \times \\
    &\qquad\qquad \left [
     b^-_{j,k\pm\frac{1}{2}}b^+_{j,k\pm\frac{1}{2}} (\Delta q)^{S(N)}_{j,k\pm1}-(b^\mp_{j,k\pm\frac{1}{2}})^2 (\Delta q)^{N(S)}_{j,k}\right] \\
    & + \frac{b^\mp_{j,k\pm\frac{1}{2}}}{(b^+_{j,k\pm\frac{1}{2}}-b^-_{j,k\pm\frac{1}{2}})\Delta y} \left[ G((\Delta q)^{S(N)}_{j,k\pm1})-G((\Delta q)^{N(S)}_{j,k}) \right],
\end{split}
\end{equation}
where 
\begin{equation} \label{C.8}
\begin{split}
    F((\Delta q)^W_{j+1,k}) &= f((\Delta q)^W_{j+1,k}+\tilde q_{j+\frac{1}{2},k})-f(\tilde q_{j+\frac{1}{2},k}), \\
    F((\Delta q)^E_{j,k}) &= f((\Delta q)^E_{j,k}+\tilde q_{j+\frac{1}{2},k})-f(\tilde q_{j+\frac{1}{2},k}), 
\end{split}
\end{equation}
and
\begin{equation} \label{C.9}
\begin{split}
    G((\Delta q)^S_{j,k+1}) &= g((\Delta q)^S_{j,k+1}+\tilde q_{j,k+\frac{1}{2}})-g(\tilde q_{j,k+\frac{1}{2}}), \\
    G((\Delta q)^N_{j,k}) &= g((\Delta q)^N_{j,k}+\tilde q_{j,k+\frac{1}{2}})-g(\tilde q_{j,k+\frac{1}{2}}).
\end{split}
\end{equation}
Here, the notation $(\Delta q)^I$ with $I=\{E,W,N,S\}$ is defined similarly to \eqref{4.10} by
\begin{equation}  \label{4.62}
\begin{split}
    & (\Delta q)^E_{j,k}:= (\Delta q)^n_{j,k} + \frac{\Delta x}{2}((\Delta q)_x)^n_{j,k},\;\;\;\;
    (\Delta q)^W_{j,k}:= (\Delta q)^n_{j,k} - \frac{\Delta x}{2}((\Delta q)_x)^n_{j,k},\\
    & (\Delta q)^N_{j,k}:= (\Delta q)^n_{j,k} + \frac{\Delta y}{2}((\Delta q)_y)^n_{j,k},\;\;\;\;
    (\Delta q)^S_{j,k}:= (\Delta q)^n_{j,k} - \frac{\Delta y}{2}((\Delta q)_y)^n_{j,k},
\end{split}
\end{equation}
with the $MC-\theta$ limiter 
\begin{equation} \label{4.63}
\begin{split}
    &((\Delta q)_x)^n_{j,k}:= \\
    & 
    \text{minmod}\left( \theta\frac{(\Delta q)^n_{j,k}-(\Delta q)^n_{j-1,k}}{\Delta x}, 
    \frac{(\Delta q)^n_{j+1,k}-(\Delta q)^n_{j-1,k}}{2\Delta x},
    \theta \frac{(\Delta q)^n_{j+1,k}-(\Delta q)^n_{j,k}}{\Delta x} \right), \\
    &((\Delta q)_x)^n_{j,k}:= \\
    & 
    \text{minmod}\left( \theta\frac{(\Delta q)^n_{j,k}-(\Delta q)^n_{j,k-1}}{\Delta y}, 
    \frac{(\Delta q)^n_{j,k+1}-(\Delta q)^n_{j,k-1}}{2\Delta y},
    \theta \frac{(\Delta q)^n_{j,k+1}-(\Delta q)^n_{j,k}}{\Delta y} \right),
\end{split}
\end{equation}
and $1\leq\theta\leq 2$. The detailed computation of \eqref{4.59} is similar to the derivation in section 3.3 of \cite{ref8}.

Next, we consider the rest terms on the right-hand-side of \eqref{4.57}. The subdomain $D_{j,k}$ can be regarded as a rectangle when $\Delta t \to 0$, up to small corners of a negligible size $O((\Delta t)^2)$. Applying \eqref{4.48} to the rest terms yields that
\begin{equation} \label{4.64}
\begin{split}
    \lim_{\Delta t \to 0}  \frac{1}{\Delta t} 
    \bigg [ & \frac{1}{\Delta x \Delta y}  |D_{j,k}|  \overline{w}^{n+1}_{j,k}  - \Delta q^n_{j,k}  \bigg ] = \\
    & \left[ \frac{a^-_{j+\frac{1}{2},k}}{\Delta x} (\Delta q)^E_{j,k} 
    + \frac{b^-_{j,k+\frac{1}{2}}}{\Delta y} (\Delta q)^N_{j,k}
    - \frac{a^+_{j-\frac{1}{2},k}}{\Delta x} (\Delta q)^W_{j,k}
    - \frac{b^+_{j,k-\frac{1}{2}}}{\Delta y} (\Delta q)^S_{j,k} \right]\\
    & - \frac{1}{\Delta x} \left[F((\Delta q)^E_{j,k}) - F((\Delta q)^W_{j,k}) \right]
    -\frac{1}{\Delta y} \left[G((\Delta q)^N_{j,k}) - G((\Delta q)^S_{j,k}) \right] \\
    & + S((\Delta q)_{j,k}) ]    
\end{split}
\end{equation}
Finally, combining \eqref{4.59} and \eqref{4.64}, the semi-discrete scheme takes the form 
\begin{equation} \label{C.13}
    \frac{d}{dt} (\Delta q)_{j,k}(t) = 
    - \frac{H^x_{j+\frac{1}{2},k}(t)-H^x_{j-\frac{1}{2},k}(t)}{\Delta x}
    - \frac{H^y_{j,k+\frac{1}{2}}(t)-H^y_{j,k-\frac{1}{2}}(t)}{\Delta y}
    + S((\Delta q)_{j,k}(t))
\end{equation}
with the numerical fluxes
\begin{equation} \label{C.14}
\begin{split}
    H^x_{j+\frac{1}{2},k}  = & \frac{a^+_{j+\frac{1}{2},k}F((\Delta q)^E_{j,k})-a^-_{j+\frac{1}{2},k}F((\Delta q)^W_{j+1,k})}{a^+_{j+\frac{1}{2},k}-a^-_{j+\frac{1}{2},k}}  \\
    &\qquad\qquad\qquad\qquad\qquad + \frac{a^+_{j+\frac{1}{2},k}a^-_{j+\frac{1}{2},k}}{a^+_{j+\frac{1}{2},k}-a^-_{j+\frac{1}{2},k}} \left[(\Delta q)^W_{j+1,k}-(\Delta q)^E_{j,k} \right], \\
    H^y_{j,k+\frac{1}{2}} = & \frac{b^+_{j,k+\frac{1}{2}}G((\Delta q)^N_{j,k})-b^-_{j,k+\frac{1}{2}}G((\Delta q)^S_{j,k+1})}{b^+_{j,k+\frac{1}{2}}-b^-_{j,k+\frac{1}{2}}} \\
    &\qquad\qquad\qquad\qquad\qquad + \frac{b^+_{j,k+\frac{1}{2}}b^-_{j,k+\frac{1}{2}}}{b^+_{j,k+\frac{1}{2}}-b^-_{j,k+\frac{1}{2}}} \left[(\Delta q)^S_{j,k+1}-(\Delta q)^N_{j,k} \right].
\end{split}
\end{equation}

\section{Maximum Principle} \label{append.D}

We begin with the explicit form of \eqref{4.69},
\begin{equation} \label{4.72}
\begin{split}
    (\Delta q)^{n+1}_{j,k} = (\Delta q)^n_{j,k} 
    & - \frac{\lambda^n}{(a^+_{j+\frac{1}{2},k}-a^-_{j+\frac{1}{2},k})} 
    \left[ a^+_{j+\frac{1}{2},k}F((\Delta q)^E_{j,k}) - a^-_{j+\frac{1}{2},k}F((\Delta q)^W_{j+1,k})  \right]  \\
    & + \frac{\lambda^n}{(a^+_{j-\frac{1}{2},k}-a^-_{j-\frac{1}{2},k})} 
    \left[ a^+_{j-\frac{1}{2},k}F((\Delta q)^E_{j-1,k}) - a^-_{j-\frac{1}{2},k}F((\Delta q)^W_{j,k})  \right] \\
    & - \frac{\lambda^n(a^+_{j+\frac{1}{2},k}a^-_{j+\frac{1}{2},k})}{a^+_{j+\frac{1}{2},k}-a^-_{j+\frac{1}{2},k}} \left[ (\Delta q)^W_{j+1,k}-(\Delta q)^E_{j,k} \right] \\
    & + \frac{\lambda^n(a^+_{j-\frac{1}{2},k}a^-_{j-\frac{1}{2},k})}{a^+_{j-\frac{1}{2},k}-a^-_{j-\frac{1}{2},k}} \left[ (\Delta q)^W_{j,k}-(\Delta q)^E_{j-1,k} \right] \\
    & - \frac{\mu^n}{(b^+_{j,k+\frac{1}{2}}-b^-_{j,k+\frac{1}{2}})} 
    \left[ b^+_{j,k+\frac{1}{2}}G((\Delta q)^N_{j,k}) - b^-_{j,k+\frac{1}{2}}G((\Delta q)^S_{j,k+1})  \right]  \\
    & + \frac{\mu^n}{(b^+_{j,k-\frac{1}{2}}-b^-_{j,k-\frac{1}{2}})} 
    \left[ b^+_{j,k-\frac{1}{2}}G((\Delta q)^N_{j,k-1}) - b^-_{j,k-\frac{1}{2}}G((\Delta q)^S_{j,k})  \right] \\
    & - \frac{\mu^n(b^+_{j,k+\frac{1}{2}}b^-_{j,k+\frac{1}{2}})}{b^+_{j,k+\frac{1}{2}}-b^-_{j,k+\frac{1}{2}}} \left[ (\Delta q)^S_{j,k+1}-(\Delta q)^N_{j,k} \right] \\
    & + \frac{\mu^n(b^+_{j,k-\frac{1}{2}}b^-_{j,k-\frac{1}{2}})}{b^+_{j,k-\frac{1}{2}}-b^-_{j,k-\frac{1}{2}}} \left[ (\Delta q)^S_{j,k}-(\Delta q)^N_{j,k-1} \right].
\end{split}
\end{equation}
All the terms on the right-hand of \eqref{4.72} are taken at the time step $t^n$.\quad By the definition \eqref{4.62}, we have the equality
\begin{equation} \label{4.73}
    (\Delta q)^n_{j,k} = \frac{(\Delta q)^E_{j,k}+(\Delta q)^W_{j,k}+(\Delta q)^N_{j,k}+(\Delta q)^S_{j,k}}{4}.
\end{equation}
Then we substitute \eqref{4.73} for $(\Delta q)^n_{j,k}$ and adjust the other terms in \eqref{4.72},
\begin{align} \label{4.74}
    (\Delta q)^{n+1}_{j,k}  
     =  &\frac{(\Delta q)^E_{j,k} +(\Delta q)^W_{j,k}+(\Delta q)^N_{j,k}+(\Delta q)^S_{j,k}}{4} \notag \\
     & + \lambda^n \left[  \frac{a^-_{j+\frac{1}{2},k}}{a^+_{j+\frac{1}{2},k}-a^-_{j+\frac{1}{2},k}} 
     \Big[ F((\Delta q)^W_{j+1,k}) - F((\Delta q)^E_{j,k}) \Big] - F((\Delta q)^E_{j,k}) \right] \notag \\
     & - \lambda^n \left[ \frac{a^+_{j-\frac{1}{2},k}}{a^+_{j-\frac{1}{2},k}-a^-_{j-\frac{1}{2},k}} 
     \Big[ F((\Delta q)^W_{j,k}) - F((\Delta q)^E_{j-1,k}) \Big] - F((\Delta q)^W_{j,k})\right] \notag \\
     & - \lambda^n \frac{(a^+_{j+\frac{1}{2},k}a^-_{j+\frac{1}{2},k})}{a^+_{j+\frac{1}{2},k}-a^-_{j+\frac{1}{2},k}} \Big[ (\Delta q)^W_{j+1,k}-(\Delta q)^E_{j,k} \Big] \notag \\
     & + \lambda^n \frac{(a^+_{j-\frac{1}{2},k}a^-_{j-\frac{1}{2},k})}{a^+_{j-\frac{1}{2},k}-a^-_{j-\frac{1}{2},k}} \Big[ (\Delta q)^W_{j,k}-(\Delta q)^E_{j-1,k} \Big] \notag \\
     & + \mu^n \left[ \frac{b^-_{j,k+\frac{1}{2}}}{b^+_{j,k+\frac{1}{2}}-b^-_{j,k+\frac{1}{2}}} 
     \Big[ G((\Delta q)^S_{j,k+1}) - G((\Delta q)^N_{j,k}) \Big]  - G((\Delta q)^N_{j,k}) \right] \notag \\
     & - \mu^n \left[ \frac{b^+_{j,k-\frac{1}{2}}}{b^+_{j,k-\frac{1}{2}}-b^-_{j,k-\frac{1}{2}}} 
     \Big[ G((\Delta q)^S_{j,k}) + G((\Delta q)^N_{j,k-1})  \Big] - G((\Delta q)^S_{j,k})   \right] \notag \\
     & - \mu^n \frac{(b^+_{j,k+\frac{1}{2}}b^-_{j,k+\frac{1}{2}})}{b^+_{j,k+\frac{1}{2}}-b^-_{j,k+\frac{1}{2}}} 
     \Big[ (\Delta q)^S_{j,k+1}-(\Delta q)^N_{j,k} \Big] \notag \\
     & + \mu^n \frac{(b^+_{j,k-\frac{1}{2}}b^-_{j,k-\frac{1}{2}})}{b^+_{j,k-\frac{1}{2}}-b^-_{j,k-\frac{1}{2}}} 
     \Big[ (\Delta q)^S_{j,k}-(\Delta q)^N_{j,k-1} \Big].
\end{align}
To simplify notations, we use the abbreviations used in \cite{ref3} 
\begin{equation} \label{4.75}
\begin{split}
     \Delta^x_{j+\frac{1}{2},k} (\Delta q)
    & := (\Delta q)^W_{j+1,k}(t^n) - (\Delta q)^E_{j,k}(t^n), \\
     \Delta^x_{j,k}F
    & :=F((\Delta q)^E_{j,k}) - F((\Delta q)^W_{j,k}) \\
     \Delta^x_{j,k}G
    & :=G((\Delta q)^N_{j,k}) - G((\Delta q)^S_{j,k}).
\end{split}
\end{equation}
Then we rewrite the \eqref{4.74} as  
\begin{align} \label{4.76}
    (\Delta q)^{n+1}_{j,k}  
     = & \frac{(\Delta q)^E_{j,k} +(\Delta q)^W_{j,k}+(\Delta q)^N_{j,k}+(\Delta q)^S_{j,k}}{4} \notag \\
     & +\lambda^n \left[ \frac{a^-_{j+\frac{1}{2},k}}{a^+_{j+\frac{1}{2},k}-a^-_{j+\frac{1}{2},k}}
     \frac{\Delta^x_{j+\frac{1}{2},k}F}{\Delta^x_{j+\frac{1}{2},k}(\Delta q)}
     \Big[ (\Delta q)^W_{j+1,k}-(\Delta q)^E_{j,k} \Big] \right. \notag \\
     & \qquad\quad - \frac{\Delta^x_{j,k}F}{\Delta^x_{j,k}(\Delta q)} \Big[(\Delta q)^E_{j,k}-(\Delta q)^W_{j,k} \Big] \notag \\
     & \left. \qquad\quad - \frac{a^+_{j-\frac{1}{2},k}}{a^+_{j-\frac{1}{2},k}-a^-_{j-\frac{1}{2},k}}
     \frac{\Delta^x_{j-\frac{1}{2},k}F}{\Delta^x_{j-\frac{1}{2},k}(\Delta q)}
     \Big[ (\Delta q)^W_{j,k}-(\Delta q)^E_{j-1,k} \Big]     \right] \notag \\
     & - \lambda^n \frac{(a^+_{j+\frac{1}{2},k}a^-_{j+\frac{1}{2},k})}{a^+_{j+\frac{1}{2},k}-a^-_{j+\frac{1}{2},k}} \Big[ (\Delta q)^W_{j+1,k}-(\Delta q)^E_{j,k} \Big] \notag \\
     & + \lambda^n \frac{(a^+_{j-\frac{1}{2},k}a^-_{j-\frac{1}{2},k})}{a^+_{j-\frac{1}{2},k}-a^-_{j-\frac{1}{2},k}} \Big[ (\Delta q)^W_{j,k}-(\Delta q)^E_{j-1,k} \Big] \notag \\
     & +\mu^n \left[ \frac{b^-_{j,k+\frac{1}{2}}}{b^+_{j,k+\frac{1}{2}}-b^-_{j,k+\frac{1}{2}}}
     \frac{\Delta^y_{j,k+\frac{1}{2}}G}{\Delta^y_{j,k+\frac{1}{2}}(\Delta q)}
     \Big[ (\Delta q)^S_{j,k+1}-(\Delta q)^N_{j,k} \Big] \right. \notag \\
     & \qquad\quad - \frac{\Delta^y_{j,k}G}{\Delta^y_{j,k}(\Delta q)} \Big[(\Delta q)^N_{j,k}-(\Delta q)^S_{j,k} \Big] \notag \\
     & \left. \qquad\quad - \frac{b^+_{j,k-\frac{1}{2}}}{b^+_{j,k-\frac{1}{2}}-b^-_{j,k-\frac{1}{2}}}
     \frac{\Delta^y_{j,k-\frac{1}{2}}G}{\Delta^y_{j,k-\frac{1}{2}}(\Delta q)}
     \Big[ (\Delta q)^S_{j,k}-(\Delta q)^N_{j,k-1} \Big]     \right] \notag \\
     & - \mu^n \frac{(b^+_{j,k+\frac{1}{2}}b^-_{j,k+\frac{1}{2}})}{b^+_{j,k+\frac{1}{2}}-b^-_{j,k+\frac{1}{2}}} 
     \Big[ (\Delta q)^S_{j,k+1}-(\Delta q)^N_{j,k} \Big] \notag \\
     & + \mu^n \frac{(b^+_{j,k-\frac{1}{2}}b^-_{j,k-\frac{1}{2}})}{b^+_{j,k-\frac{1}{2}}-b^-_{j,k-\frac{1}{2}}} 
     \Big[ (\Delta q)^S_{j,k}-(\Delta q)^N_{j,k-1} \Big]
\end{align}
Collecting the coefficients of $(\Delta q)^W_{j+1,k}, (\Delta q)^{E(W)}_{j,k}, (\Delta q)^E_{j-1,k}$, and $(\Delta q)^S_{j,k+1}$, $(\Delta q)^{N(S)}_{j,k}$, $(\Delta q)^N_{j,k-1}$, 
\begin{equation} \label{4.77}
\begin{split}
    (\Delta q)^{n+1}_{j,k} =  
    & \lambda^n \left(
    \frac{a^-_{j+\frac{1}{2},k}}{a^+_{j+\frac{1}{2},k}-a^-_{j+\frac{1}{2},k}}
    \frac{\Delta^x_{j+\frac{1}{2},k}F}{\Delta^x_{j+\frac{1}{2},k}(\Delta q)}
    - \frac{a^+_{j+\frac{1}{2},k}a^-_{j+\frac{1}{2},k}}{a^+_{j+\frac{1}{2},k}-a^-_{j+\frac{1}{2},k}}
    \right) (\Delta q)^W_{j+1,k} \\
    & + \left[ \frac{1}{4} - \lambda^n \left(
    \frac{a^-_{j+\frac{1}{2},k}}{a^+_{j+\frac{1}{2},k}-a^-_{j+\frac{1}{2},k}}
    \frac{\Delta^x_{j+\frac{1}{2},k}F}{\Delta^x_{j+\frac{1}{2},k}(\Delta q)}
    +  \frac{\Delta^x_{j,k}F}{\Delta^x_{j,k}(\Delta q)} \right. \right. \\
    &\qquad\qquad\qquad\qquad\qquad\qquad\qquad\qquad\qquad \left.\left.- \frac{a^+_{j+\frac{1}{2},k}a^-_{j+\frac{1}{2},k}}{a^+_{j+\frac{1}{2},k}-a^-_{j+\frac{1}{2},k}}    
    \right) \right] (\Delta q)^E_{j,k} \\
    & + \left[\frac{1}{4} + \lambda^n \left(
    \frac{\Delta^x_{j,k}F}{\Delta^x_{j,k}(\Delta q)}
    - \frac{a^+_{j-\frac{1}{2},k}}{a^+_{j-\frac{1}{2},k}-a^-_{j-\frac{1}{2},k}}
    \frac{\Delta^x_{j-\frac{1}{2},k}F}{\Delta^x_{j-\frac{1}{2},k}(\Delta q)} \right.\right.\\
    &\qquad\qquad\qquad\qquad\qquad\qquad\qquad\qquad\qquad \left.\left. + \frac{a^+_{j-\frac{1}{2},k}a^-_{j-\frac{1}{2},k}}{a^+_{j-\frac{1}{2},k}-a^-_{j-\frac{1}{2},k}} 
    \right)    \right] (\Delta q)^W_{j,k} \\
    & +\lambda^n \left(
    \frac{a^+_{j-\frac{1}{2},k}}{a^+_{j-\frac{1}{2},k}-a^-_{j-\frac{1}{2},k}}
    \frac{\Delta^x_{j-\frac{1}{2},k}F}{\Delta^x_{j-\frac{1}{2},k}(\Delta q)}
    - \frac{a^+_{j-\frac{1}{2},k}a^-_{j-\frac{1}{2},k}}{a^+_{j-\frac{1}{2},k}-a^-_{j-\frac{1}{2},k}}
    \right) (\Delta q)^E_{j-1,k} \\   
    & + \mu^n \left(
    \frac{b^-_{j,k+\frac{1}{2}}}{b^+_{j,k+\frac{1}{2}}-b^-_{j,k+\frac{1}{2}}}
    \frac{\Delta^y_{j,k+\frac{1}{2}}G}{\Delta^y_{j,k+\frac{1}{2}}(\Delta q)}
    - \frac{b^+_{j,k+\frac{1}{2}}b^-_{j,k+\frac{1}{2}}}{b^+_{j,k+\frac{1}{2}}-b^-_{j,k+\frac{1}{2}}}
    \right) (\Delta q)^S_{j,k+1} \\
    & + \left[ \frac{1}{4} - \mu^n \left(
    \frac{b^-_{j,k+\frac{1}{2}}}{b^+_{j,k+\frac{1}{2}}-b^-_{j,k+\frac{1}{2}}}
    \frac{\Delta^y_{j,k+\frac{1}{2}}G}{\Delta^y_{j,k+\frac{1}{2}}(\Delta q)}
    +  \frac{\Delta^y_{j,k}G}{\Delta^y_{j,k}(\Delta q)} \right. \right. \\
    &\qquad\qquad\qquad\qquad\qquad\qquad\qquad\qquad\qquad \left.\left.- \frac{b^+_{j,k+\frac{1}{2}}b^-_{j,k+\frac{1}{2}}}{b^+_{j,k+\frac{1}{2}}-b^-_{j,k+\frac{1}{2}}}    
    \right) \right] (\Delta q)^N_{j,k} \\
    & + \left[\frac{1}{4} + \mu^n \left(
    \frac{\Delta^y_{j,k}G}{\Delta^y_{j,k}(\Delta q)}
    - \frac{b^+_{j,k-\frac{1}{2}}}{b^+_{j,k-\frac{1}{2}}-b^-_{j,k-\frac{1}{2}}}
    \frac{\Delta^y_{j,k-\frac{1}{2}}G}{\Delta^y_{j,k-\frac{1}{2}}(\Delta q)} \right.\right.\\
    &\qquad\qquad\qquad\qquad\qquad\qquad\qquad\qquad\qquad \left.\left. + \frac{b^+_{j,k-\frac{1}{2}}b^-_{j,k-\frac{1}{2}}}{b^+_{j,k-\frac{1}{2}}-b^-_{j,k-\frac{1}{2}}} 
    \right)    \right] (\Delta q)^S_{j,k} \\
    & +\mu^n \left(
    \frac{b^+_{j,k-\frac{1}{2}}}{b^+_{j,k-\frac{1}{2}}-b^-_{j,k-\frac{1}{2}}}
    \frac{\Delta^y_{j,k-\frac{1}{2}}G}{\Delta^y_{j,k-\frac{1}{2}}(\Delta q)}
    - \frac{b^+_{j,k-\frac{1}{2}}b^-_{j,k-\frac{1}{2}}}{b^+_{j,k-\frac{1}{2}}-b^-_{j,k-\frac{1}{2}}}
    \right) (\Delta q)^N_{j,k-1}.
\end{split}
\end{equation}
Next, we discuss the coefficients. Due to the fact that $F'(\Delta q) = f'(q)$ (see lemma \ref{lemma1}) and the fact that $a^\pm_{j+\frac{1}{2},k}$ is the maximal speed on its direction, which is determined from the flux derivatives, (consult the definition \eqref{4.9}), we obtain the following inequalities:
\begin{equation} \label{4.78}
\begin{split}
    a^+_{j+\frac{1}{2},k}\geq 0 \qquad 
    & \text{and} \qquad
    |\frac{\Delta^x_{j+\frac{1}{2},k}F}{\Delta^x_{j+\frac{1}{2},k}(\Delta q)}| \leq a^+_{j+\frac{1}{2},k},\\    
    a^-_{j+\frac{1}{2},k} \leq 0, \qquad
    & \text{and} \qquad
    |\frac{\Delta^x_{j+\frac{1}{2},k}F}{\Delta^x_{j+\frac{1}{2},k}(\Delta q)}| \leq - a^-_{j+\frac{1}{2},k}.
\end{split}
\end{equation}
Hence, the  rearranged coefficients of $(\Delta q)^W_{j+1,k}$ and $(\Delta q)^E_{j-1,k}$ are non-negative,
\begin{equation} \label{4.79}
\begin{split}
    & \underbrace{\lambda^n \frac{(-a^-_{j+\frac{1}{2},k})}{a^+_{j+\frac{1}{2},k}-a^-_{j+\frac{1}{2},k}}}_{\geq0} \;
    \underbrace{\Bigg[
     a^+_{j+\frac{1}{2},k} 
     - \frac{\Delta^x_{j+\frac{1}{2},k}F}{\Delta^x_{j+\frac{1}{2},k}(\Delta q)}
    \Bigg]}_{\geq0} \geq 0 \\ 
    & \underbrace{\lambda^n \frac{a^+_{j-\frac{1}{2},k}}{a^+_{j-\frac{1}{2},k}-a^-_{j-\frac{1}{2},k}}}_{\geq0} \;
    \underbrace{\Bigg[ 
     \frac{\Delta^x_{j-\frac{1}{2},k}F}{\Delta^x_{j-\frac{1}{2},k}(\Delta q)}
     - a^-_{j-\frac{1}{2},k}
     \Bigg]}_{\geq0} \geq 0.
\end{split}
\end{equation}
The coefficients of $(\Delta q)^E_{j,k}$ and $(\Delta q)^W_{j,k}$ are also non-negative due to the CFL assumption \eqref{4.70}.\quad By the assumption \eqref{4.70}, we have
\begin{equation}
    \lambda^n a^+_{j+\frac{1}{2},k}\leq \frac{1}{8} \quad \text{and} \quad
    -\frac{1}{8} \leq \lambda^n a^-_{j+\frac{1}{2},k}.   
\end{equation}
The above inequalities imply
\begin{equation}
\begin{split}
    &\lambda^n a^+_{j+\frac{1}{2},k} - \lambda^n a^-_{j+\frac{1}{2},k}  \leq \frac{1}{4} \\
    \Rightarrow \quad
    & 4 \leq \frac{1}{\lambda^n a^+_{j+\frac{1}{2},k}- \lambda^n a^-_{j+\frac{1}{2},k}} \\
    \Rightarrow \quad & -\frac{1}{2}\leq\frac{\lambda^n a^-_{j+\frac{1}{2},k}}{\lambda^n a^+_{j+\frac{1}{2},k}- \lambda^n a^-_{j+\frac{1}{2},k}} \\
    \Rightarrow \quad & \frac{-\lambda^n a^-_{j+\frac{1}{2},k}}{\lambda^n a^+_{j+\frac{1}{2},k}- \lambda^n a^-_{j+\frac{1}{2},k}} \leq \frac{1}{2}
\end{split}
\end{equation}
Hence, the coefficient of $(\Delta q)^E_{j,k}$ is non-negative because of  
\begin{equation} \label{4.80}
\begin{split}
    & \frac{1}{4}+\underbrace{ \frac{(-a^-_{j+\frac{1}{2},k})}{a^+_{j+\frac{1}{2},k}-a^-_{j+\frac{1}{2},k}}}_{\leq \frac{1}{2}}
    \Bigg( \underbrace{\lambda^n\frac{\Delta^x_{j+\frac{1}{2},k}F}{\Delta^x_{j+\frac{1}{2},k}(\Delta q)}}_{\geq-\frac{1}{8}}
    \underbrace{- \lambda^n a^+_{j+\frac{1}{2},k}}_{\geq -\frac{1}{8}} \Bigg)
    \underbrace{-\lambda^n \frac{\Delta^x_{j,k}F}{\Delta^x_{j,k}(\Delta q)}}_{\geq -\frac{1}{8}} \geq0. 
\end{split}
\end{equation}
By the same way, the coefficient of $(\Delta q)^W_{j,k}$ can be proved to be non-negative, and the other four coefficients as well.\quad Since all the coefficients are non-negative and the sum of the coefficients are equal to 1, the combination on the right-hand-side of \eqref{4.77} is a convex combination.\quad Hence,
\begin{equation*}
    (\Delta q)^{n+1}_{j,k}\leq 
    \max \Big((\Delta q)^W_{j+1,k}, (\Delta q)^{E(W)}_{j,k}, (\Delta q)^E_{j-1,k}, (\Delta q)^S_{j,k+1}, (\Delta q)^{N(S)}_{j,k}, (\Delta q)^N_{j,k-1} \Big).
\end{equation*}
Because of the definition of the intermediate values $(\Delta q)^{E,W,S,N}$ and the choice of the derivative in \eqref{4.62}, these intermediate values satisfy the local maximum principle, (consult Theorem 1 of \cite{ref11}),
\begin{equation*}
\begin{split}
    \max_{j,k} \Big((\Delta q)^W_{j+1,k}, (\Delta q)^{E(W)}_{j,k}, (\Delta q)^E_{j-1,k}, (\Delta q)^S_{j,k+1}, (\Delta q)^{N(S)}_{j,k},& (\Delta q)^N_{j,k-1} \Big) \\
    & \leq \max_{j,k}((\Delta q)^n_{j,k}).
\end{split}
\end{equation*}
Then the maximum principle: $\max_{j,k}(\Delta q)^{n+1}_{j,k}\leq \max_{j,k}(\Delta q)^n_{j,k}$ holds. 
\hfill $\blacksquare$

\end{document}